\documentclass[12pt]{article}

\usepackage[latin1]{inputenc}
\usepackage{amsmath,amssymb}
\usepackage{graphics, graphicx,psfrag}
\usepackage{amsthm}
\usepackage{epsfig}
\usepackage{amscd}
\usepackage{verbatim} 
\usepackage[english,frenchb]{babel}

\newcommand{\R}{R}
\newcommand{\Q}{Q}  
\newcommand{\T}{T}
\newcommand{\F}{F} 

\newcommand{\B}{\mathcal B} 
\newcommand{\N}{\mathcal N} 

\newcommand{\G}{\Gamma}

\newcommand{\I}{\mathrm{I}}
\newcommand{\IIi}{\mathrm{II_1}}
\newcommand{\II}{\mathrm{II}}
\newcommand{\III}{\mathrm{III}}

\newcommand{\CI}{{\mathrm{\bf C}}}
\newcommand{\RI}{{\mathrm{\bf R}}}

\newcommand{\ZI}{{\mathrm{\bf Z}}}
\newcommand{\NI}{{\mathrm{\bf N}}}

\newcommand{\SI}{{\boldsymbol S}}

\newcommand{\SL}{\mathrm {SL}}

\newcommand{\Int}{\mathrm{Int}}

\newcommand{\Aut}{\mathrm{Aut}}

\newcommand{\Out}{\mathrm{Out}}
\newcommand{\Def}{\mathrm{Def}}

\newcommand{\norm}[1]{||#1||}
\newcommand{\abs}[1]{|#1|}

\newcommand{\graph}{\mathrm{graph}}

\newcommand{\impl}{\Longrightarrow}
\newcommand{\ssi}{\Longleftrightarrow}

\newcommand{\del}{\partial}

\newcommand{\id}{\mathrm{Id}}
\renewcommand{\Im}{\mathop{\mathrm{Im}}}

\newcommand{\h}{{\mathfrak{h}}} 

\newtheorem*{dfn}{\indent\bf D\'efinition}
\newtheorem{thm}{\indent\bf Th\'eor\`eme}
\newtheorem{lem}[thm]{\indent\bf Lemme}

\newtheorem{prop}[thm]{\indent\bf Proposition}
\newtheorem{cor}[thm]{\indent\bf Corollaire}

\renewcommand{\qed}{_\blacksquare}
\newenvironment{dem}{\noindent D\'emonstration. }{{\hfill {$\qed$}}\vspace{.25cm} }

\makeatletter
\renewcommand\thepart{\@Roman\c@part.}
\renewcommand\thesection{\@arabic\c@section.}
\renewcommand\thesubsection{\@arabic\c@subsection.}
\renewcommand\section{\@startsection {section}{2}{\z@}%
    {0ex \@plus -1ex  \@minus -.2ex}%
    {2.3ex \@plus .2ex}%
    {\reset@font \center \sc }}
\renewcommand\subsection{\@startsection {subsection}{3}{\parindent}%
    {.1ex \@plus .5ex \@minus .2ex}%
    {-.8em \@plus .2ex}%
    { \reset@font\small \bfseries}}
\makeatother

\usepackage{diagrams}

\begin{document}

\centerline{\bf{\uppercase{Sur les espaces mesur\'es singuliers}}}

\medskip

\centerline{{\bf I - \'Etude m\'etrique-mesur\'ee}\footnote{
The author was partially supported by a JSPS Fellowship for European Researchers.}}
\vspace{1cm}
\centerline{\small {Mika\"el \sc Pichot}}
\vspace{.4cm}
\centerline{\small{Novembre 2004}}
\vspace{.2cm}
\centerline{---}
\vspace{.3cm}

\selectlanguage{english}
\begin{abstract}
Recall Jones-Schmidt theorem that an ergodic measured equivalence relation is strongly ergodic if and only if it has no nontrivial amenable quotient. We give two new characterizations of strong ergodicity, in terms of metric-measured spaces. 
The first one identifies strong ergodicity with the concentration property as defined, in  this  (foliation) setting, by Gromov \cite{Gromov00_SQ}. 
The second one characterize the existence of nontrivial amenable quotients in terms of ``F\o lner sequences" in graphings naturally associated to (the leaf space of) the equivalence relation.  

We also present a formalization of the concept of quasi-periodicity, based on (noncommutative) measure theory. The ``singular measured spaces" appearing in the title refer to the leaf spaces of measured equivalence relations.
\end{abstract}
\selectlanguage{french}

\tableofcontents


\newpage

\section{Introduction}

L'importance fondamentale des structures transverses de feuilletages a \'et\'e mise en avant au milieu des ann\'ees 1950 par les travaux d'Andr\'e Haefliger (\cite{Heafliger58} par exemple, voir l'expos\'e historique \cite{Heafliger03}).  Ces structures, d\'efinies sur les vari\'et\'es transverses et invariantes par holonomie, refl\`etent des propri\'et\'es de \textit{l'espace des feuilles} du feuilletage consid\'er\'e.\\

Il est bien connu que l'espace des feuilles en lui-m\^eme est, le plus souvent, \og singulier\fg. Ainsi en est-il des exemples les plus simples, dont les feuilletages lin\'eaires du tore $\RI^2/\ZI^2$ de dimension 2, par droites de pente irrationnelle (feuilletages de Kronecker). N\'eanmoins,  une analyse non triviale de ces espaces reste possible. Elle a \'et\'e initi\'ee \`a la fin des ann\'ees 1970 par A. Connes dans \cite{Connes79}. De fa\c con g\'en\'erale, le concept d'\textit{espace singulier} traduit l'existence sur de nombreux espaces quotients \textit{a priori} non standard,  en particulier sur l'espace des feuilles d'un feuilletage ou d'une lamination, de structures canoniques non triviales et intrins\`eques. Ces structures, de m\^eme que pour les espaces classiques, peuvent \^etre de natures diverses, e.g. mesur\'ee, topologique, diff\'erentielle, m\'etrique. \\

Cette \'etude concerne la \textit{th\'eorie de la mesure des espaces singuliers}. Elle est divis\'ee en deux parties.  La premi\`ere partie (le pr\'esent article)    traite de  la notion d'\textit{ergodicit\'e forte}, en relation avec les \textit{structures m\'etriques-mesur\'ees}  naturellement associ\'ees aux espaces singuliers. La seconde partie   \cite{Pichot04_II} traite de  la  \textit{propri\'et\'e T de Kazhdan}  pour les espaces mesur\'es singuliers.

\bigskip
\bigskip

\centerline{---}

\bigskip

Les \textit{relations d'\'equivalence} interviennent par nature m\^eme dans l'\'etude des espaces quotients. Elles en constituent, par d\'efinition, des  \textit{d\'e\-sin\-gu\-la\-ri\-sa\-tions}. D\'esingulariser un espace singulier en un espace usuel a pour avantage imm\'ediat d'en permettre l'\'etude \`a l'aide d'outils math\'ematiques standard. D'une certaine mesure, ce point de vue rapproche les espaces singuliers des vari\'et\'es : une d\'e\-sin\-gu\-la\-ri\-sa\-tion est l'analogue d'une carte en g\'eom\'etrie diff\'erentielle, son r\^ole est de (sur)param\'etrer convenablement l'espace quotient. Une structure singuli\`ere sur un ensemble consiste, ainsi, en la donn\'ee d'un \og syst\`eme de d\'esingularisations compatibles\fg\ de cet ensemble. La th\'eorie de la mesure des espaces singuliers repose alors sur la notion de \textit{relation d'\'equivalence mesur\'ee}, fond\'ee dans \cite{FeldmanMoore77}. 

\bigskip

Soit $X$ un espace bor\'elien standard. Une relation d'\'equivalence \textit{\`a classes d\'enom\-brables}  sur $X$ est \textit{bor\'elienne} si son graphe $\R\subset X\times X$  est r\'eunion d'une famille d\'enombrable d'isomorphismes partiels bor\'eliens de $X$. Lorsque $X$ est muni d'une mesure de probabilit\'e sans atome $\mu$ et que ces isomorphismes partiels  \textit{pr\'eservent la classe de $\mu$}, on dit que $\R$ est une \textit{relation d'\'equivalence mesur\'ee} (cf. \textsection\ref{rel}). Par exemple, une action $\alpha$ d'un groupe d\'enombrable par automorphismes bor\'eliens de $(X,\mu)$ d\'efinit une relation d'\'equivalence bor\'elienne sur $X$ (la partition en orbites), et cette relation est mesur\'ee si $\mu$ est \textit{quasi-invariante} par $\alpha$. Plusieurs travaux importants de diff\'erents auteurs, notamment Adams, Furman, Gaboriau, Hjorth, Kechris, Monod, Popa,  Shalom,  ont r\'ecemment contribu\'e  \`a un d\'eveloppement significatif de la th\'eorie.

\bigskip

Deux relations d'\'equivalence mesur\'ees $R_1$ sur $(X_1,\mu_1)$  et $R_2$ sur $(X_2,\mu_2)$ d\'esingu\-larisent un m\^eme espace mesur\'e singulier si et seulement si elles sont \textit{stablement isomorphes}, au sens o\`u il existe un isomorphisme bor\'elien non singulier $\rho : X_1' \to X_2'$ entre deux parties bor\'eliennes non n\'egligeables $X_1'\subset X_1$  et $X_2'\subset X_2$ tel que $x\sim_{\R_1} y \ssi \rho(x)\sim_{\R_2} \rho(y)$ pour $x,y \in X_1'$.

Nous dirons d'une propri\'et\'e de relation d'\'equivalence mesur\'ee, invariante par isomorphisme stable, qu'elle est une \textit{propri\'et\'e de l'espace mesur\'e singulier} des orbites de cette relation.

L'\textit{ergodicit\'e} et l'\textit{ergodicit\'e forte} sont des propri\'et\'es d'espaces mesur\'es singuliers.
 L'ergodicit\'e est une notion dynamique classique.  Une relation d'\'equivalence mesur\'ee sur $X$ est \textit{ergodique} si tout bor\'elien \textit{satur\'e}  est n\'egligeable ou de compl\'ementaire n\'egligeable, i.e. si $X$ ne contient pas de bor\'eliens invariants non triviaux.  
  L'ergodicit\'e forte, \'egalement de nature dynamique, a \'et\'e introduite par Connes-Weiss \cite{ConnesWeiss80} et Schmidt \cite{Schmidt80}. Une action bor\'elienne $\alpha$ d'un groupe d\'enombrable $\G$ quasi-pr\'eservant une mesure de probabilit\'e $\mu$ sur $X$ est  \textit{fortement ergodique} si elle ne poss\`ede pas de parties bor\'eliennes asymptotiquement invariantes non triviales. Rappelons que, suivant \cite{ConnesKrieger77,Schmidt80,ConnesWeiss80,Schmidt81}, une suite $A_n\subset X$ de parties bor\'eliennes est dite \textit{asymptotiquement invariante} sous l'action $\alpha$ si 
\[
\mu(\alpha(\gamma)A_n\Delta A_n) \to_n 0
\]
pour tout $\gamma \in \Gamma$, et qu'elle est dite \textit{non triviale} s'il existe $\delta >0$ tel que 
\[
\delta \leqslant \mu(A_n) \leqslant 1- \delta.
\]
On v\'erifie que l'existence de suites asymptotiquement invariantes non triviales ne d\'epend effectivement que de l'espace singulier des orbites de l'action $\alpha$.

\bigskip

La notion d'ergodicit\'e forte a permis \`a Connes-Weiss \cite{ConnesWeiss80} de caract\'eriser les groupes de Kazhdan d\'enombrables par leurs actions ergodiques \textit{pr\'eservant une mesure de probabilit\'e} (que l'on appellera ergodique de type $\IIi$).  Schmidt  \cite{Schmidt81} et Losert-Rindler \cite{LosertRindler81} ont ensuite caract\'eris\'e les groupes moyennables,  dans le m\^eme esprit. Les r\'esultats  sont les suivants : 
\begin{itemize}
\item[1.]  K. Schmidt \cite{Schmidt80}  a observ\'e  que toute action ergodique de type $\IIi$ d'un groupe de Kazhdan est fortement ergodique.
\item[2.] Un groupe qui poss\`ede pas la propri\'et\'e T poss\`ede au moins une action ergodique de type $\IIi$ non fortement ergodique \cite{ConnesWeiss80}.
\item[3.]  Un groupe est moyennable si et seulement s'il ne poss\`ede pas d'action ergodique de type $\IIi$ qui soit fortement ergodique \cite{LosertRindler81,Schmidt81}. 
\end{itemize}

\bigskip

Citons \'egalement un th\'eor\`eme particuli\`erement remarquable, le \textit{th\'eor\`eme de Jones-Schmidt}, montrant qu'une relation d'\'equivalence mesur\'ee  (\`a classes d\'enombrables) ergodique est fortement ergodique si et seulement si elle ne poss\`ede pas de quotient moyennable non trivial (cf. \cite{JonesSchmidt87}).\\

L'article r\'ecent G. Hjorth et A. Kechris \cite[App. 1]{HjorthKechris03} contient une pr\'esentation  d\'etaill\'ee de la notion d'ergodicit\'e forte pour les relations d'\'equivalence mesur\'ees.

\bigskip
\bigskip

\centerline{---}

\bigskip

Avant d'\'enoncer nos r\'esultats, d\'ecrivons la  nature des structures m\'etriques-mesur\'ees  associ\'ees aux espaces singuliers.

\bigskip

Au niveau intuitif, la situation est la suivante. Consid\'erons un complexe simplicial  $Y$, au sens usuel du terme, non compact, mais satisfaisant \`a certaines propri\'et\'es de \textit{quasi-p\'eriodicit\'e}. On peut alors  l\'egitimement esp\'erer associer, \`a chaque partie bor\'elienne quasi-p\'eriodique  de $Y$, un nombre r\'eel positif repr\'esentant le \textit{covolume} de cette partie, de valeur proportionelle \`a sa \og densit\'e\fg\ dans $Y$. Dans le cas p\'eriodique, i.e. lorsque $Y$ est muni d'une action cocompacte d'un groupe discret, la donn\'ee d'une mesure sur l'espace quotient d\'etermine une telle densit\'e. La co-existence d'un covolume et d'une m\'etrique (de la m\'etrique simpliciale par exemple) d\'efinit alors sur $Y$ une structure \textit{d'espace m\'etrique-mesur\'e}. De  fa\c con similaire, \textit{les espaces m\'etriques-mesur\'es associ\'es aux espaces singuliers sont des espaces m\'etriques quasi-p\'eriodiques, munis d'un covolume}.

\bigskip

Soit $M$ une vari\'et\'e et $F$ un feuilletage \textit{minimal} sur $M$ (i.e. toutes les feuilles sont denses). Il est relativement clair que la fonction indicatrice d'un voisinage ouvert d'un point de $M$ est, en restriction \`a chaque feuille de $F$, une fonction \og quasi-p\'eriodique\fg\ d\'efinie sur cette feuille. Si $F$ est \textit{ergodique}, la m\^eme observation est valable, non seulement pour les voisinages ouverts, mais aussi pour la fonction indicatrice de toute partie bor\'elienne non n\'egligeable ; celle-ci d\'efinit une fonction quasi-p\'eriodique en restriction \`a une feuille g\'en\'erique de $F$. Ces observations sont  bien connues. Elles peuvent \^etre illustr\'ees plus pr\'ecis\'ement, par exemple, par les th\'eor\`emes de Ghys \cite{Ghys95} et Cantwell-Conlon \cite{CantwellConlon98} sur la topologie des feuilles d'une feuilletage de dimension 2. Il est \'egalement int\'eressant de noter, r\'eciproquement, qu'un espace m\'etrique poss\'edant certaines propri\'et\'es de quasi-p\'eriodicit\'e peut parfois \^etre plong\'e dans une lamination minimale, ou ergodique, de fa\c con interessante  ; pour un exemple concret d'une telle construction, nous renvoyons le lecteur \`a l'\'etude des quasi-cristaux telle qu'expos\'ee dans \cite{Bellissard} par exemple. 

\bigskip

Le point de vue que nous adoptons dans cet article est de consid\'erer que chaque espace singulier, par exemple l'espace $\Q = M/F$ des feuilles de $F$, est un concept (une notion) de quasi-p\'eriodicit\'e ---  nous nous int\'eressons uniquement ici aux \textit{aspects mesur\'es} du concept de quasi-p\'eriodicit\'e et de sa formalisation (i.e. au cas \textit{ergodique}). 

\bigskip

Consid\'erons l'exemple le plus simple, celui des \textit{graphes quasi-p\'eriodiques associ\'es \`a un espace mesur\'e singulier $\Q$} (que l'on appellera  \textit{graphes $\Q$-p\'eriodiques}) : 

\textbullet)  Soit $\R$ une relation d'\'equivalence mesur\'ee sur $(X,\mu)$. Soit $K\subset \R$ une partie bor\'elienne sym\'etrique de $\R$. On dit que $K$ est un \textit{graphage} de $\R$ si pour tous points \'equivalents $x,y\in X$, il existe un nombre fini $x_0=x,x_1,\ldots,x_n=y$ de points de $X$ tels que $(x_i,x_{i+1})\in K$. En suspendant $K$ au-dessus de $X$, i.e. en attachant pour tout $(x,y)\in K$ une ar\^ete $\simeq [0,1]$ entre $x$ et $y$,  on obtient alors une \og lamination\fg\ $\Sigma_K$ dont les feuilles sont des graphes, et dont $X$ est une transversale totale. Soulignons que les orbites de $\R$ sont connexes dans $\Sigma_K$.  Connes, Feldman, et Weiss \cite{CFW81} ont pour la premi\`ere fois consid\'er\'e  les parties bor\'eliennes $K\subset \R$ comme des familles mesurables de graphes sur les orbites de $\R$ ; la notion de graphage au sens pr\'esent\'e ici (i.e. lorsque les orbites sont connexes)  a \'et\'e introduite par Levitt  (cf. \cite{Gaboriau99}). (voir \textsection\ref{rel}) 

\textbullet)  \`A une relation d'\'equivalence mesur\'ee $\R$ sur $X$ et un graphage  $K$ de $\R$, on associe  le \textit{graphe quasi-p\'eriodique} $\tilde \Sigma_K$  dont, par d\'efinition,  les sommets sont les points de $\R\subset X\times X$, et les ar\^etes sont les couples $((x,y),(x,z))$ de points de $\R$ tels que $(y,z)\in K$. Soit $\Q$ un espace singulier. On dit que $\tilde \Sigma_K$ est \textit{un graphe $\Q$-p\'eriodique} lorsque $\Q$ est l'espace de ses feuilles (i.e. $\Q=X/\R$).

La relation entre $\tilde \Sigma_K$ et $\Sigma_K$ est simple : $\Sigma_K$ est le quotient du graphe  quasi-p\'eriodique $\tilde \Sigma_K$ par la relation de quasi-p\'eriodicit\'e ($\Sigma_K =\tilde \Sigma_K/\R$).  Ces d\'efinitions s'\'etendent bien s\^ur en dimension sup\'erieure. Les complexes simpliciaux $\Q$-p\'eriodiques ont \'et\'e introduits par D. Gaboriau dans \cite{Gaboriau02} sous le nom de $\R$-complexe simplicial  (cf. \textsection\ref{QP}). Plus g\'en\'eralement, de nombreuses cat\'egories d'espaces m\'etriques s\'eparables  contiennent de fa\c con naturelle  des espaces quasi-p\'eriodiques en un sens analogue (en particulier les vari\'et\'es riemanniennes). 

\bigskip

La \textit{structure m\'etrique-mesur\'ee} sur $\Sigma_K$ est par d\'efinition donn\'ee par la m\'etrique simpliciale $d$ sur chaque feuille et par la mesure transverse $\mu$ d\'efinie sur les bor\'eliens  de $X$. Cette structure se rev\`ele \`a $\tilde \Sigma_K$ qui est l'exemple fondamental d'espace m\'etrique-mesur\'e quasi-p\'eriodique associ\'e \`a $\Q$. Le support du covolume ainsi d\'efini sur  $\tilde \Sigma_K$ est constitu\'e de bor\'eliens \textit{quasi-p\'eriodiques} (de sommets, i.e. inclus dans $\R$).

\bigskip

Notons que la m\'etrique simpliciale sur $\tilde \Sigma_K$ peut prendre la valeur $+\infty$ (deux points sont \`a distance infinie si et seulement s'ils ne sont pas sur une m\^eme feuille) : les complexes simpliciaux quasi-p\'eriodiques associ\'es aux espaces mesur\'es singuliers sont des espaces m\'etriques-mesur\'es au sens o\`u l'a d\'efini M. Gromov dans \cite{Gromov00_SQ}  (dans le cadre des feuilletages, o\`u la m\'etrique est longitudinale et la mesure, la mesure de Lebesgue sur la vari\'et\'e).

\bigskip

Ce paragraphe peut \^etre r\'esum\'e par les diagrammes suivants.

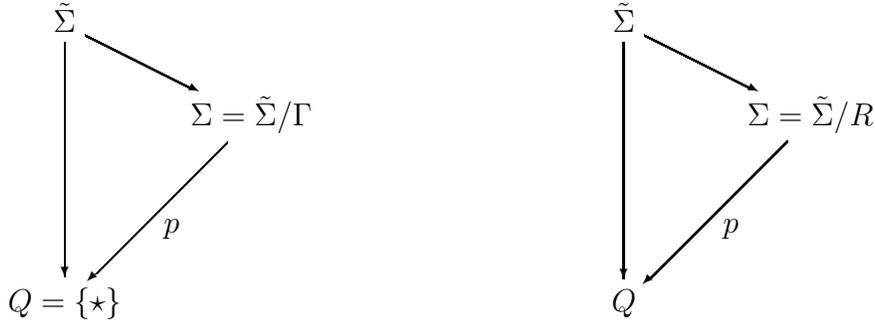
\begin{diagram}[nohug,notextflow]
\tilde \Sigma &             &               &&&    & \tilde \Sigma      &      &   \\
\dTo     &\rdTo(2,1) & \Sigma=\tilde \Sigma/\G  &&&  & \dTo  & \rdTo(2,1) &  \Sigma=\tilde\Sigma/\R  &\\
        &    \ldTo_p     &            &&&  &  & \ldTo_p    &                         &\\
\Q=\{\star\} &         &              &&& & \Q &             &             \\
\end{diagram}

\bigskip
\noindent Sur le digramme de gauche, $\tilde \Sigma$ est un complexe simplicial usuel muni d'une action libre cocompacte de d'un groupe discret $\G$. Il s'agit du cas p\'eriodique.  Sur le diagramme de droite   $\tilde \Sigma$ est un complexe simplicial quasi-p\'eriodique au sens ci-dessus muni d'une action  libre d'une relation d'\'equivalence $\R$. Il s'agit du cas quasi-p\'eriodique.  L'hypoth\`ese de cocompacit\'e est remplac\'ee dans le cas quasi-p\'eriodique par la d\'efinition suivante.

\bigskip

\begin{dfn}
On dit qu'un complexe simplicial $\Q$-p\'eriodique est \textnormal{uniform\'ement localement fini (u.l.f.)} si le nombre de simplexes attach\'es en chacun de ses sommets est uniform\'ement fini.
\end{dfn}

 \bigskip

\centerline{---}

\bigskip

L'\'etude des espaces  m\'etriques-mesur\'es est en premier lieu concern\'ee par la question de la concentration de la mesure. Notre premier r\'esultat \'etablit un lien direct entre ergodicit\'e forte et concentration.

\bigskip

La plus simple des notions de quasi-p\'eriodicit\'e (non p\'eriodique) obtenue \`a l'aide des espaces mesur\'es singuliers est \textit{l'hyperfinitude}. Consid\'erons par exemple une \textit{droite quasi-p\'eriodique} $\Sigma$, i.e. $\Sigma = \Sigma_K$ est associ\'e au graphe $K \subset X\times X$ d'un automorphisme partiel $\varphi$ de $(X,\mu)$ agissant essentiellement librement (automorphisme ap\'eriodique) en pr\'eservant la mesure $\mu$. L'espace singulier $\Q=X/\langle \varphi\rangle$ est hyperfini (ce qui signifie que la relation d'\'equivalence $\R_\varphi$ des orbites de $\varphi$ est r\'eunion \textit{croissante} de sous-relations \`a \textit{orbites finies}).  On peut alors v\'erifier la propri\'et\'e suivante. 
\'Etant donn\'es un nombre fini de parties bor\'eliennes quasi-p\'eriodiques de $\Sigma$ de covolume suffisament petit, on peut \'ecarter (pour la m\'etrique simpliciale) chacune de ces parties les unes des autres sans toutefois en modifier les covolumes respectifs.
En effet, le \textit{lemme de Rokhlin} montre qu'il existe pour tout $n$ une famille finie $\Omega_1,\ldots \Omega_n$ de bor\'eliens disjoints de $X$ tels que $\varphi(\Omega_i)=\Omega_{i+1}$ et $\mu(X\backslash \amalg \Omega_i) \leqslant \varepsilon_n$ o\`u $\varepsilon_n \to_n 0$. Il est alors \'el\'ementaire de construire deux (ou plusieurs) suites $(A_n)$ et $(B_n)$ de mesure (convergeant vers) 1/4 (ou $c>0$ assez petit), form\'ees de r\'eunions de certaines des parties $\Omega_i$,  et telles que la distance entre $A_n$ et $B_n$ converge vers $+\infty$.

\bigskip

\begin{dfn}[\cite{Gromov00_SQ}] On dit qu'un espace m\'etrique-mesur\'e  $(\Sigma, d, \mu)$ est \textnormal{concentr\'e} s'il existe une fonction $c : ]0,1]^2 \to \RI_+$, telle que pour tous bor\'eliens non n\'egligeables $A,B \subset \Sigma$, on a
\[
d(A,B) \leqslant c(\mu(A),\mu(B)).
\]
\end{dfn}

\bigskip

Nous dirons d'une propri\'et\'e d'espace singulier ergodique qu'elle est un \textit{param\`etre de quasi-p\'eriodicit\'e}.  Le r\'esultat qui suit (cf. \textsection\ref{concentr})  montre que la concentration des graphes $\Q$-p\'eriodiques  un param\`etre de quasi-p\'eriodicit\'e (l'hypoth\`ese \og de type fini\fg\ figurant dans ce th\'eor\`eme est comparable \`a celle utilis\'ee en th\'eorie des groupes et refl\`ete l'existence de suffisamment de complexes simpliciaux $Q$-p\'eriodiques u.l.f.).

\bigskip

\begin{thm}
Soit $\Q$ un espace singulier ergodique de type fini. Les conditions suivantes sont \'equivalentes,
\begin{itemize}
\item [i.] Il existe un graphe $\Q$-p\'eriodique u.l.f. concentr\'e,
\item [ii.] Tout graphe $\Q$-p\'eriodique u.l.f. est concentr\'e.
\end{itemize}
\end{thm}

\bigskip

On dira dans ce cas que l'espace singulier $\Q$ est concentr\'e.

\bigskip

\textit{Observation.}   Ce th\'eor\`eme s'\'etend \`a d'autres espaces m\'etriques s\'eparables que les graphes. Rappelons que la structure m\'etrique-mesur\'ee consid\'er\'ee pour les graphes est donn\'ee par le couple (distance simpliciale,  covolume) ; le support du covolume est de dimension 0 au sens o\`u il n\'eglige les bor\'eliens (non quasi-p\'eriodiques ou) non port\'es par les sommets.  Il est n\'eanmoins clair que, si ces graphes sont concentr\'es, alors  les complexes simpliciaux $\Q$-p\'eriodiques u.l.f., munis par exemple d'un covolume port\'e sur les simplexes de dimension $i$, sont \'egalement concentr\'es ; il en r\'esulte aussi, en adaptant les techniques de triangulation des vari\'et\'es au cas quasi-p\'eriodique (cf. \cite{HeitschLazarov91}), que les vari\'et\'es $\Q$-p\'eriodiques sont \'egalement concentr\'ees, o\`u l'on remplace l'hypoth\`ese u.l.f. par l'hypoth\`ese \og \`a g\'eom\'etrie born\'ee\fg.  (cf \cite{Pichot05})

\bigskip

Nous montrons au paragraphe \ref{concentr} le th\'eor\`eme suivant.
 
\bigskip

\begin{thm}\label{thm1}
Un espace singulier ergodique de type fini est concentr\'e si et seulement s'il est fortement ergodique.
\end{thm}

Le point de vue de la concentration de la mesure pour les relations d'\'equivalence, et notamment une version \textit{quantitative} de la propri\'et\'e de concentration des espaces quasi-p\'eriodiques associ\'es \`a un espace singulier, sont susceptibles de donner de nouveaux invariants pour les relations fortement ergodiques. Rappelons que, dans le cas p\'eriodique, la propri\'et\'e de concentration est un ph\'enom\`ene trivial au niveau qualitatif (tous les espaces sont concentr\'es car compacts), mais que des ph\'enom\`enes remarquables apparaissent au niveau quantitatif (e.g. concentration normal ou exponentielle). Il semble donc int\'eressant de mener une \og syst\'ematique \fg\ de la concentration pour les relations d'\'equivalence mesur\'ees. Cette \'etude est en projet.

\bigskip
\bigskip

\centerline{---}

\bigskip

La \textit{classification des relations d'\'equivalence moyennables} (mesur\'ees \`a classes d\'enom\-brables) a \'et\'e achev\'ee en 1981 suite \`a la d\'emonstration par Connes, Feldman et Weiss du th\'eor\`eme suivant. \textit{Toute relation d'\'equivalence mesur\'ee moyennable peut \^etre engendr\'ee par une seule transformation de l'espace.} Les auteurs montrent l'\'equivalence de plusieurs d\'efinitions du concept de moyennabilit\'e pour les relations d'\'equivalences mesur\'ees, et obtiennent notamment une caract\'erisation des relations moyennables en terme de suites de F\o lner pr\'esentes dans les structures de graphes mesurables sur les orbites de la relation (cf. \cite{CFW81}). Nous d\'efinissons au paragraphe \ref{isop}, \'etant donn\'e un graphe quasi-p\'eriodique $\tilde \Sigma$, la notion de \textit{suites de F\o lner \'evanescentes} dans $\tilde \Sigma$.  Il s'agit d'une reformulation g\'eom\'etrique de la notion dynamique de \textit{I-suites} consid\'er\'ee par Schmidt dans \cite{Schmidt81} pour des actions de groupes d\'enombrables pr\'eservant une mesure de probabilit\'e.

\begin{dfn} Soit $\tilde\Sigma$ un graphe $\Q$-p\'eriodique u.l.f. muni de sa structure m\'etrique-mesur\'ee $( \tilde d, \tilde \mu)=$(m\'etrique simpliciale, 0-covolume) d\'ecrite ci-dessus.  \'Etant donn\'ee une partie bor\'elienne $A\subset \tilde \Sigma$, on note $\del_{K} A$ l'ensemble des points de $\tilde \Sigma \backslash A$ \`a distance $\leqslant 1$ de $A$ dans $\tilde \Sigma$. On dit qu'une suite $(A_n)$ de parties bor\'eliennes non n\'egligeables de $\tilde \Sigma$ est une   \textnormal{suite de F\o lner \'evanescente}  s'il existe  une suite $(\varepsilon_n)$ de nombres r\'eels convergeant vers 0  telles que  
\[
\tilde \mu(A_n) \to 0 \qquad \mathrm{et} \qquad \tilde \mu(\del_{K} A_n) \leqslant \varepsilon_n \tilde \mu(A_n).
\]
\end{dfn}

\bigskip

Inspir\'es par \cite{CFW81}, nous montrons le th\'eor\`eme suivant.

\bigskip

\begin{thm}
Un espace singulier ergodique de type fini $\Q$ poss\`ede un quotient moyennable si et seulement si tout graphe $\Q$-p\'eriodique uniform\'ement localement fini poss\`ede des suites de F\o lner \'evanescentes.
\end{thm}

\bigskip
Ce th\'eor\`eme \'egalement s'\'etend \`a d'autres espaces m\'etriques-mesur\'es que les graphes.

\bigskip
\bigskip

\centerline{---}
\bigskip

Je remercie Damien Gaboriau pour son aide constante au cours de l'\'elaboration de ce travail.

Je dois \'egalement beaucoup \`a \'Etienne Ghys, ainsi qu'aux excellentes conditions de travail dont on b\'en\'eficie au sein de l'UMPA.

Je remercie Yann Ollivier pour sa lecture critique du manuscrit. Je remercie Damien Gaboriau pour ses lectures, du manuscrit et des nombreuses versions ant\'erieures, qui auront permis d'am\'eliorer consid\'erablement ce texte.

\vspace{1cm}
\section{Relations d'\'equivalence mesur\'ees}\label{rel}

Soit $X$ un espace bor\'elien standard. Rappelons qu'il s'agit d'un espace polonais (topologique s\'eparable admettant une m\'etrique compl\`ete) muni de sa structure bor\'elienne. On dit qu'une relation d'\'equivalence $\R$ sur $X$ est \textit{bor\'elienne} si son graphe $ \R\subset X\times X$ est une partie bor\'elienne.

\bigskip

Les relations d'\'equivalence \textit{\`a classes d\'enombrables} jouent un r\^ole privil\'egi\'e dans la th\'eorie. Un r\'esultat bien connu de Feldman-Moore montre qu'une telle relation est bor\'elienne si et seulement si on peut munir \textit{mesurablement} chacune de ses orbites d'une structure de graphe complet. Plus pr\'ecis\'ement, il existe une partition d\'enombrable 
\[
\R = \amalg_{i\in \NI} \graph(\varphi_i)\subset X\times X
\]
de toute relation bor\'elienne $\R$, \`a classe d\'enombrables, en graphes d'isomorphismes partiels $\varphi_i : A_i \xrightarrow{\simeq} B_i$ entre deux bor\'eliens $A_i$ et $B_i$ de $X$ (cf. \cite{FeldmanMoore77}) ; on  \textit{construit} ainsi entre deux points \'equivalents quelconques $x$ et $y$ de $X$ une unique ar\^ete orient\'ee $x \xrightarrow{\varphi_i}y$ \cite{Gaboriau99,Gaboriau02} (ce qui, d'un point de vue alg\'ebrique, revient \`a postuler l'existence d'un espace classifiant pour $\R$).

\bigskip

\textit{Exemples.} i. Une action $\alpha : \G \to \Aut(X)$ d'un groupe d\'enombrable $\G$, par isomorphismes bor\'eliens,  d\'efinit sur $X$ une relation d'\'equivalence bor\'elienne $\R_\alpha$ \`a classes d\'enombrables, donn\'ee par la partition de $X$ en les orbites de $\alpha$. Le r\'esultat de Feldman-Moore ci-dessus peut s'interpr\^eter en disant que toute relation d'\'equivalence \`a classes d\'enombrables est la partition en orbites d'une action mesurable de groupe discret (en choisissant une partition de cette relation par des isomorphismes partiels d'ordre 2 \'etendus \`a $X$). 

ii. Un feuilletage sur une vari\'et\'e la partitionne en feuilles et d\'efinit ainsi une relation d'\'equivalence bor\'elienne. En restreignant cette relation \`a une transversale $T$, on obtient une relation d'\'equivalence \`a classes d\'enombrables (dont le graphe $\R \subset T\times T$ est partitionn\'e par des applications d'holonomie). 

iii. Une source importante de relations d'\'equivalence provient, par leur nature m\^eme, des probl\`emes de classifications. On se contentera ici d'\'evoquer un exemple, l'espace des groupes de type fini, \'etudi\'e dans \cite{Champetier00}. L'espace $X$ consid\'er\'e est l'espace topologique compact des groupes marqu\'es (l'espace des quotients d'un groupe libre), sur lequel on \'etudie  (par exemple) la relation d'isomorphisme.  Il s'agit d'une relation d'\'equivalence bor\'elienne \`a classes d\'enombrables.

\bigskip

Soit $\R$ une relation d'\'equivalence bor\'elienne sur $X$. On appelle \textit{isomorphisme partiel (int\'erieur) de $\R$} un isomorphisme partiel $\varphi : A \to B$ entre deux bor\'eliens de $X$ tel que $\varphi (x) \sim x$  pour tout $x\in A$. L'ensemble des isomorphismes partiels de $\R$ se note $[[\R]]$. 
Le groupe des automorphismes int\'erieurs d'une relation d'\'equivalence se note $[\R]$, ou $\Int(\R)$, et s'appelle le \textit{groupe plein}. Il est ainsi constitu\'e des isomorphismes $X\to X$ dont le graphe est inclus dans $\R \subset X\times X$. Il s'agit bien s\^ur d'un invariant d'isomorphisme au sens suivant.

\medskip

\begin{dfn}
Deux relations d'\'equivalence bor\'eliennes $\R$ et $\R'$ sur $X$  et $X'$ sont \textnormal{isomorphes} s'il existe un isomorphisme bor\'elien $\rho : X \to X'$  tel que $x\sim y$ si et seulement si $\rho(x)\sim' \rho(y)$ pour tout $x,y\in X$. On dit que $\R$ et $\R'$ sont \textnormal{stablement isomorphes} s'il existe deux parties bor\'eliennes $\Omega\subset X$ et $\Omega' \subset X'$ rencontrant respectivement toutes les classes de $\R$ et $\R'$, telles que les relations restreintes $\R_{|\Omega}$ et $\R'_{|\Omega'}$ soient isomorphes.
\end{dfn}

\medskip

\textit{Exemple.} Un feuilletage et sa restriction \`a une transversale totale d\'efinissent deux relations d'\'equivalence stablement isomorphes. 

\bigskip

\bigskip

La th\'eorie des relations d'\'equivalence dans le cadre bor\'elien est principalement d\'evelopp\'ee en logique (cf. \cite{JKL} par exemple). Dans la suite de ce texte (et dans \cite{Pichot04_II}), nous supposerons toujours la pr\'esence additionnelle d'une \textit{mesure quasi-invariante}, dans la tradition de \cite{FeldmanMoore77} (et de Murray-von  Neumann) --- on identifie alors deux relations d'\'equivalence ayant \textit{presque} les m\^emes orbites.\\

 Plus pr\'ecis\'ement, soit $(X,\mu)$ un espace de probabilit\'e, i.e. un espace bor\'elien standard muni d'une mesure bor\'elienne de probabilit\'e sans atome. On dit qu'une relation d'\'equivalence bor\'elienne \textit{\`a classes d\'enombrables} sur $X$ est une \textit{relation d'\'equivalence mesur\'ee} si la mesure $\mu$ est \textit{quasi-invariante} au sens o\`u tout bor\'elien n\'egligeable a un satur\'e n\'egligeable (notons que le satur\'e d'un bor\'elien, i.e. la r\'eunion des classes intersectant ce bor\'elien,  est bor\'elien). Le th\'eor\`eme c\'el\`ebre suivant, par exemple, constitue une question encore ouverte dans le cadre bor\'elien (avec une hypoth\`ese convenable d'irr\'eductibilit\'e rempla\c  cant l'ergodicit\'e).

\medskip

\begin{thm}[Connes-Feldman-Weiss \cite{CFW81}] Soit $\R$ une relation d'\'equivalence ergodique moyennable sur $(X,\mu)$. Il existe un isomorphisme  bor\'elien $T$ de l'espace $X$ pr\'eservant la classe de $\mu$ tel que $x\sim_\R y \ssi \ y = T^n(x)$ pour presque tous $x,y\in X$.
\end{thm} 

\medskip

Rappelons qu'une relation d'\'equivalence mesur\'ee est dite \textit{ergodique} si les parties bor\'eliennes invariantes (i.e. satur\'ees) sont n\'egligeables ou co-n\'egligeables. La notion de moyennabilit\'e intervenant dans le th\'eor\`eme a \'et\'e d\'efinie par Zimmer. Toute action d'un groupe discret moyennable pr\'eservant la classe de $\mu$ d\'efinie une relation d'\'equivalence mesur\'ee moyennable au sens de Zimmer.  Rappelons \'egalement qu'une relation d'\'equivalence est engendr\'ee par un seul automorphisme de $X$ si et seulement si elle est \textit{hyperfinie}, i.e. si elle s'\'ecrit comme r\'eunion d\'enombrable croissante de relations d'\'equivalence mesur\'ees \`a classes finies (Dye). Le th\'eor\`eme de Connes-Feldman-Weiss est la conclusion d'une s\'erie de travaux, dont ceux d'Ornstein-Weiss sur la g\'en\'eralisation du lemme de Rokhlin aux groupes moyennables et le th\'eor\`eme selon lequel une action  ergodique d'un groupe moyennable,  pr\'eservant une mesure de probabilit\'e, est hyperfinie \cite{OrnsteinWeiss80}.

\bigskip

\begin{dfn}
Deux relations d'\'equivalence mesur\'ees $\R$ et $\R'$ sur $(X,\mu)$  et $(X',\mu')$ sont \textnormal{isomorphes} (resp. \textnormal{stablement isomorphes}) s'il existe deux bor\'eliens $\Omega \subset X$ et $\Omega' \subset X'$ de mesure totale (resp. dont les satur\'es sont de mesure totale) et un isomorphisme de relations d'\'equivalence bor\'eliennes entre $\R_{|\Omega}$ et $\R_{|\Omega'}$, qui est non singulier au sens o\`u il envoie la classe de $\mu$ sur la classe de $\mu'$.
\end{dfn}

\bigskip
Nous renvoyons \`a \cite{Mackey66,Connes79,Ramsay82} pour l'extension de ces notions aux relations d'\'equivalence \`a classes non n\'ecessairement d\'enombrables.

\bigskip

Une relation d'\'equivalence mesur\'ee ergodique peut \^etre de type $\II$ ou de type $\III$, selon qu'il existe ou non une mesure $\sigma$-finie invariante dans la classe de $\mu$. Une mesure quasi-invariante $\sigma$-finie $\mu$  est dite \textit{invariante pour $\R$} si pour une partition $\R=\amalg_i \graph(\varphi_i)$ en graphes d'isomorphismes partiels, on a $\mu(\varphi_i(\Omega))=\mu(\Omega)$ pour tout $\Omega$ inclus dans le domaine de $\varphi_i$. V\'erifier que cette d\'efinition est ind\'ependante de la partition choisie constitue un exercice typique de la th\'eorie g\'eom\'etrique des relations d'\'equivalence mesur\'ees, nous renvoyons \`a \cite{Gaboriau99} pour de nombreuses illustrations de cette technique (d\'ecoupage des domaines). Lorsqu'il existe une mesure de probabilit\'e invariante dans la classe de $\mu$, on dit que $\R$ est de type $\IIi$. Ces d\'efinitions s'\'etendent aux relations non ergodiques.

\bigskip

\textit{Remarque.} La notion d'isomorphisme d\'ecrite ci-dessus (\og l'\'equivalence orbitale \fg), ainsi que la r\'epartition en types,  ont \'et\'e introduites par Murray et von Neumann au cours de leurs travaux sur les alg\`ebres d'op\'erateurs (1936-1943). L'alg\`ebre de von Neumann associ\'ee \`a une action libre d'un groupe d\'enombrable sur $(X,\mu)$ ne d\'epend que (de la relation d'\'equivalence mesur\'ee form\'ee) des orbites de cette action.

\bigskip
\bigskip

Concluons cette section par des faits standard.\\

Soit $\R$ une relation d'\'equivalence mesur\'ee sur un espace de probabilit\'e $(X,\mu)$. La mesure $\mu$ s'\'etend canoniquement \`a $\R$ en une mesure $\h$ d\'efinie par
\[
\h(K) = \int_X \#K^x d\mu(x)
\]
o\`u $K\subset \R$ est une partie bor\'elienne et $K^x = \{(x,y)\in K\cap \R\}$ (mesure de d\'ecompte horizontal). L'image $\h^{-1}$ de $\h$ par l'inversion $^{-1} : \R \to \R$ d\'efinie par  $(x,y)^{-1}=(y,x)$ est \'equivalente \`a $\h$ et on obtient un homomorphisme bor\'elien $\delta$ de $\R$ (vu comme groupo\"\i de mesurable)  dans $]0,\infty[$, i.e. une application mesurable v\'erifiant
\[
\delta(x,z) = \delta(x,y)\delta(y,z)
\]
pour tous $x,y,z$ \'equivalents, en posant
\[
d\h(x,y) = \delta(x,y)d\h^{-1}(x,y)
\] 
(d\'eriv\'ee au sens de Radon-Nikodym). Notons que $\h^{-1}$ est la mesure de d\'ecompte vertical d\'efinie par $\h^{-1}(K)=\int_X \#K_y d\mu(y)$ o\`u $K_y=\{(x,y)\in K\cap\R\}$. (cf. \cite{FeldmanMoore77})\\

Nous dirons qu'une partie bor\'elienne $K\subset \R$ est \textit{sym\'etrique} si $K=K^{-1}$. Une partie bor\'elienne sym\'etrique $K\subset \R$ d\'efinit une distance $d_K : \R \to [0,\infty]$  sur les orbites de $\R$ associant \`a $(x,y)\in \R$ le plus petit entier $n$ pour lequel il existe une suite $x_0=x$, $x_1,\ldots x_{n-1}$, $x_n=y$ telle que $(x_i,x_{i+1})\in K$. On dit que $K$ est un \textit{graphage de $\R$} si $d(x,y)<\infty$ pour presque tout $(x,y)\in \R$ (on supposera toujours dans la suite qu'un graphage est une partie \textit{sym\'etrique} de $\R$, i.e. qu'il est non orient\'e). Cela revient \`a dire que $\R = \cup_n K^n$ \`a un n\'egligeable pr\`es, o\`u $K^n$ est l'ensemble des couples $(x,y)\in \R$ tels que $d(x,y)\leqslant n$, ou encore que presque toutes les orbites sont connexes pour la structure simpliciale obtenue en \og collant\fg\ une ar\^ete entre deux points $x$ et $y$ de $X$ si et seulement si $(x,y)\in K$. Un graphage peut \^etre \'etiquet\'e par une famille d\'enombrable de lettres $\Phi$,  i.e. on peut choisir une famille d\'enombrable $\Phi$ d'isomorphismes partiels de $\R$, de sorte que $K = \cup_{\varphi \in \Phi}  \graph(\varphi)$, o\`u l'\'etiquettage est bijectif si cette  r\'eunion est disjointe (on peut bien s\^ur supposer alors que $\Phi$ est sym\'etrique au sens o\`u $\varphi \in \Phi \ssi \varphi^{-1} \in \Phi$). On appellera \'egalement $\Phi$ un \textit{graphage} de $\R$ (cette d\'efinition a \'et\'e introduite par Levitt en relation avec la notion de \textit{co\^ut} pour les relations d'\'equivalence mesur\'ees de type $\IIi$, cf. \cite{Gaboriau99}). Si les \'el\'ements de $\Phi$ partitionnent $K$, la structure simpliciale associ\'ee \`a $K$ coincide avec la structure obtenue par le proc\'ed\'e standard de suspension de $\Phi$ au-dessus de $X$. (cf. \cite{FeldmanMoore77,Gaboriau99,Gaboriau02})\\

Nous dirons qu'une partie bor\'elienne sym\'etrique $K\subset \R$ est u.l.f. (uniform\'ement localement finie) si $\# K^x$ et $\#K_y$ sont uniform\'ement finis sur $X$, et qu'elle est u.l.b. (uniform\'ement localement born\'ee), relativement \`a $\mu$, si elle est u.l.f. et si $\abs \delta_K = \sup_{(x,y)\in K} \abs{ \delta(x,y)}$ est fini. Toute partie u.l.f. $K\subset \R$ peut \^etre partitionn\'ee en un nombre fini d'isomorphismes partiels de $\R$ ; toute partie u.l.b. $K \subset \R$ peut \^etre partitionn\'e en un nombre fini  d'isomorphismes partiels $\Phi =(\varphi_1,\ldots,\varphi_n)$ de $\R$ tels que les fonctions $d({\varphi_i}_*\mu)/d\mu$ soient uniform\'ement born\'ees. Nous dirons qu'une relation d'\'equivalence mesur\'ee est \textit{de type fini} si elle poss\`ede un graphage u.l.f., i.e. si elle peut \^etre engendr\'ee par un nombre fini d'isomorphismes partiels. (cf. \cite{FeldmanMoore77,CFW81} et notamment \cite[Lemma 3]{CFW81})

\vspace{1cm}
\section{Espaces singuliers}\label{sing}

Soit $\Q$ un ensemble.

\begin{dfn}
On appelle \textnormal{d\'e\-sin\-gu\-la\-ri\-sa\-tion} (bor\'elienne) de $\Q$ la donn\'ee d'un espace bor\'elien standard $X$ et d'une application surjective \textnormal{d\'efinissable} $p:X\to \Q$.
\end{dfn}

\medskip

Nous dirons qu'une application $p : X \to \Q$ est \textit{d\'efinissable} si
\[
\R_p=\{(x,y) \in X\times X \ |\ p(x)=p(y)\}
\]
est une partie bor\'elienne de $X\times X$ (i.e. $\R_p$ est une relation d'\'equivalence bor\'elienne).

\bigskip

\textit{Exemple.} Soit $M$ une vari\'et\'e et $\F$ un feuilletage de $M$. L'application 
\[
p : M \to M/\F
\]
\textit{d\'efinie} par $x\mapsto \ell$ o\`u $\ell$ est l'unique feuille contenant $x$, est une d\'e\-sin\-gu\-la\-ri\-sa\-tion de l'espace $M/\F$ des feuilles de $F$.

\bigskip

\begin{dfn}
On dit que deux d\'e\-sin\-gu\-la\-ri\-sa\-tions $p : X \to \Q$ et $p':X'\to \Q$ sont \textnormal{\'equivalentes} (au sens bor\'elien) s'il existe deux applications bor\'eliennes $\varphi : X \to X'$ et $\varphi' : X'\to X$ telles que $p'\varphi = p$ et $p\varphi' = p'$.
\end{dfn}

On v\'erifie imm\'ediatement qu'on obtient ainsi une relation d'\'equivalence sur les d\'e\-sin\-gu\-la\-ri\-sa\-tions.

\bigskip

\textit{Exemple.} Consid\'erons une lamination $L$ sur un espace topologique $X$ et notons $p :X \to X/L$ la d\'e\-sin\-gu\-la\-ri\-sa\-tion naturelle.   Soit $T\subset X$ une transversale totale de $L$. Notons  $p'=p_{|T} : T \to X/L$ la d\'e\-sin\-gu\-la\-ri\-sa\-tion associ\'ee. Alors $p$ et $p'$ sont \'equivalente. En effet $p=p'r$ o\`u $r$ est une r\'etraction mesurable $X \to T$, obtenue par exemple en fixant une famille mesurable de m\'etriques le long des feuilles et en associant (mesurablement) \`a $x\in X$ l'un des points de $T$ le plus proche de $x$.

\bigskip

\begin{dfn}
On appelle \textnormal{structure singuli\`ere} (bor\'elienne) sur $\Q$ la donn\'ee d'une classe d'\'equi\-valence de d\'e\-sin\-gu\-la\-ri\-sa\-tions. Un \textnormal{espace  singulier} (bor\'elien) est un ensemble muni d'une structure singuli\`ere.
\end{dfn}

\medskip

En pratique, l'espace singulier ainsi qu'une ou plusieurs de ses d\'e\-sin\-gu\-la\-ri\-sa\-tions apparaissent souvent de fa\c con naturelle. Citons simplement ici,

\begin{itemize}
\item l'espace des groupes de type fini (cf. \cite{Champetier00}),

\item l'espace des immeubles de type $\tilde A_2$ (cf. \cite{Pichot04_II,BarrePichot04_I,BarrePichot04_II}).
\end{itemize}
De nombreux exemples suppl\'ementaires figurent dans \cite{Connes95,Connes04}.

\bigskip
\bigskip

Soit $\Q$ un espace singulier.

\medskip

On supposera toujours que $\Q$ admet une d\'e\-sin\-gu\-la\-ri\-sa\-tion discr\`ete au sens suivant.

\begin{dfn} 
On dit qu'une d\'e\-sin\-gu\-la\-ri\-sa\-tion $p :X \to \Q$ de $\Q$ est \textnormal{discr\`ete} si les fibres de $p$ sont d\'enombrables.
\end{dfn}

(On renvoie \`a \cite{FHM78,Ramsay82} pour des r\'esultats g\'en\'eraux concernant l'existence de d\'e\-sin\-gu\-la\-ri\-sa\-tions discr\`etes \og presque s\^urement surjectives\fg.)

\bigskip

\begin{lem} \label{SOE} Les relations d'\'equivalence bor\'eliennes $\R_p$ et $\R_{p'}$ associ\'ees \`a deux d\'e\-sin\-gu\-la\-ri\-sa\-tions discr\`etes (\'equivalentes) $p : X \to \Q$ et $p' : X' \to Q$ de $\Q$ sont stablement isomorphes.
\end{lem}

\begin{dem} Soit $\varphi : X \to X'$ une application bor\'elienne telle que $p'\varphi = p$. En particulier $x\sim_p y$ si et seulement si $\varphi(x)\sim_{p'} \varphi(y)$. L'image $X_1'\subset X'$ de $\varphi$ est une partie bor\'elienne de $X'$ et on peut choisir une section bor\'elienne $s : X_1' \to X$ de $\varphi$ ($\varphi$ \'etant \`a fibres d\'enombrables). Il est clair que $\varphi : X_1 \to X_1'$ r\'ealise un isomorphisme entre les restrictions de $\R_p$ et $\R_{p'}$ respectivement \`a l'image $X_1$ de $s$ et \`a $X_1'$. Notons que $X_1$ rencontre toutes les classes de $\R_p$, et comme $p$ est surjective et v\'erifie $p'\varphi=p$, il en est de m\^eme de $X_1'$ relativement \`a $R_{p'}$. En particulier $R_p$ et $R_{p'}$ sont stablement orbitalement \'equivalentes.
\end{dem}

\medskip

\textit{Remarque.} Ce lemme reprend le fait bien connu que les restrictions d'une relation d'\'equivalence \`a deux bor\'eliens rencontrant toutes les orbites sont stablement isomorphes. Ici $\R_p$ et $\R_{p'}$ sont les restrictions de $\R_q$ \`a $X$ et $X'$, o\`u $\R_q$ est associ\'ee \`a la d\'e\-sin\-gu\-la\-ri\-sa\-tion discr\`ete $q = p \amalg p' : X\amalg X' \to \Q$. Observons que l'isomorphisme partiel construit ici est int\'erieur, en ce sens qu'il fixe l'espace singulier $\Q$.

\bigskip

\begin{dfn} Soit $\Q$ et $\Q'$ deux espaces singuliers. On appelle \textnormal{application d\'efinissable de $\Q$ vers $\Q'$} une application $\rho : \Q \to \Q'$ telle qu'il existe une application bor\'elienne $\overline \rho : X \to X'$, o\`u  $p : X\to \Q$ et $p' : X'\to \Q'$ sont deux d\'e\-sin\-gu\-la\-ri\-sa\-tions de $\Q$ et $\Q'$, v\'erifiant $\rho p = p'\overline\rho$. On dira que $\overline \rho$ d\'esingularise $\rho$.
\end{dfn}

\textit{Exemple.} Un automorphisme ext\'erieur d'une relation d'\'equivalence bor\'elienne induit une application bijective (bi-)d\'efinissable de l'espace singulier associ\'e.

\bigskip

\begin{lem}
Soit $\rho : \Q\to \Q'$ une application d\'efinissable. Il existe  une application d\'esingularisante $\overline \rho : X \to X'$ de $\rho$ entre $X$ et $X'$, o\`u $p : X \to \Q$ et $p' : X' \to \Q'$ sont deux d\'e\-sin\-gu\-la\-ri\-sa\-tions discr\`etes donn\'ees de $\Q$ de $\Q'$.
\end{lem}

\begin{dem}
Ceci r\'esulte imm\'ediatement du diagramme commutatif
\begin{diagram}[nohug,notextflow]
   X    &\rTo^{\varphi}    & Y  &   \rTo^{\tilde \rho}  & Y'& \rTo{\varphi'} & X'    \\
        &\rdTo_{p_X}      & \dTo^{p}   &                 &\dTo{p'}   &\ldTo_{p_{X' }} & \\
        &                & \Q  &   \rTo^{\rho}          & \Q'&             &      \\
\end{diagram}
o\`u $X\to \Q$ et $X'\to \Q'$ sont des d\'e\-sin\-gu\-la\-ri\-sa\-tions discr\`etes respectivement de $\Q$ et $\Q'$, et $\tilde \rho : Y \to Y'$ une application d\'esingularisante de $\rho$. Il suffit en effet de poser $\overline \rho = \varphi'\tilde \rho\varphi : X \to X'$.
\end{dem}

\bigskip
On note $\Def(\Q)$ l'ensemble des bijections d\'efinissables de $\Q$.

\begin{cor}
$\Def(\Q)$  est un groupe (groupe des automorphismes de $\Q$).
\end{cor}

\begin{dem}
Observons tout d'abord que si $\rho$ est bijective, on peut choisir une application d\'esingularisante $\overline \rho : X\to X'$ entre deux d\'e\-sin\-gu\-la\-ri\-sa\-tions discr\`etes, qui soit un isomorphisme de relations d'\'equivalence bor\'eliennes (cf. la preuve du lemme \ref{SOE}). Par suite $\Def(\Q)$ est stable par inversion. Par ailleurs si $\overline \rho : X_1\to X_2$ et $\overline \rho': X_1'\to X_2'$ sont deux d\'e\-sin\-gu\-la\-ri\-sa\-tions bijectives de $\rho$ et $\rho'$, alors il existe un bor\'elien $A\subset X_2$ et un bor\'elien $B\subset X_1'$, et un isomorphisme $\varphi : A \to B$ entre les relations ${\R_{p_2}}_{|A}$ et ${\R_{p_1'}}_{|B}$. Alors $\overline \rho' \varphi \overline \rho : \overline \rho^{-1}(A) \to \overline \rho'(B)$ est une d\'e\-sin\-gu\-la\-ri\-sa\-tion de $\rho'\rho$, et $\Def(\Q)$ est stable par produit.
\end{dem}

\bigskip

\textit{Mesure transverse sur $\Q$.} On dira que deux d\'e\-sin\-gu\-la\-ri\-sa\-tions $p:X\to \Q$ et $p':X' \to \Q$ de $\Q$ sont \textit{conjugu\'ees} s'il existe une bijection bor\'elienne $\varphi$ de $X$  sur $X'$ telle que $p'\varphi=p$. On note $\T$ la famille des ensembles bor\'eliens $X$ lorsque $X\to \Q$  parcourt les d\'e\-sin\-gu\-la\-ri\-sa\-tions discr\`etes de $\Q$. Il est clair que $\T$ est stable par r\'eunion disjointe (d\'enombrable). On appelle \textit{mesure transverse (invariante) sur $\Q$} la donn\'ee d'une application
\[
\Lambda : \B(\T) \to [0,\infty]
\]
d\'efinie sur les bor\'eliens d'\'el\'ements de $\T$ et satisfaisant aux propri\'et\'es suivantes :
\begin{itemize}
\item $\sigma$-additivit\'e, i.e. $\Lambda(\amalg \Omega_i) = \sum \Lambda(\Omega_i)$ pour toute partie bor\'elienne $\Omega_i \subset X_i$, o\`u $(p_i :X_i \to \Q)_i$ est une famille (au plus d\'enombrable) de d\'e\-sin\-gu\-la\-ri\-sa\-tions.
\item invariance, i.e. $\Lambda(X) =\Lambda(X')$ si $X$ et $X'$ sont deux d\'e\-sin\-gu\-la\-ri\-sa\-tions conjugu\'ees.
\end{itemize}

\bigskip

\textit{Parties n\'egligeables de $\Q$.} On appelle bor\'elien de $\Q$ la projection d'un bor\'elien par une d\'esingularisation discr\`ete $X\to \Q$.  La tribu ${\cal B}(\Q)$ obtenue sur $\Q$ co\"\i ncide donc avec la tribu des bor\'eliens satur\'es de $X$ et ne d\'epend pas de la d\'esingularisation discr\`ete choisie. La donn\'ee d'une mesure invariante $\Lambda$ sur $\Q$ permet de d\'efinir sans ambig\"uit\'e la notion de partie n\'egligeable $N\subset \Q$.  On notera $\N\subset {\cal B}(\Q)$ la famille des bor\'eliens n\'egligeables de $\Q$ relatifs \`a $\Lambda$. On appelera \textit{espace mesur\'e singulier} la donn\'ee d'un espace singulier $\Q$ et d'une famille de bor\'eliens n\'egligeables $\N$ associ\'ee \`a une mesure transverse invariante $\Lambda$.

\bigskip

\begin{dfn}
On dit que deux espaces mesur\'es singuliers $(\Q_1,\N_1)$ et $(Q_2,\N_2)$ sont \textnormal{isomorphes} s'il existe deux parties n\'egligeables $N_1\subset \Q_1$ et $N_2\subset Q_2$ telles que les espaces mesur\'es singuliers $\Q_1\backslash N_1$ et $\Q_2\backslash N_2$ sont \textnormal{strictement isomorphes} au sens o\`u il existe une bijection d\'efinissable non singuli\`ere $\rho$ entre les deux (i.e. $N$ est n\'egligeable si et seulement si $\rho(N)$ est n\'egligeable).
\end{dfn}

Ainsi deux d\'e\-sin\-gu\-la\-ri\-sa\-tions discr\`etes d'espaces mesur\'es singuliers isomorphes sont stablement isomorphes en tant que relations d'\'equivalence mesur\'ees. Nous dirons d'une propri\'et\'e de relation d'\'equivalence mesur\'ee, invariante par isomorphisme stable, qu'elle est une  \textit{propri\'et\'e de l'espace singulier des orbites de cette relation}.

\bigskip

L'ensemble des applications d\'efinissables non singuli\`eres de $\Q$, bijectives en restriction au compl\'emetaire d'une partie n\'egligeable, forme un groupe (en effet la composition $\varphi_2\circ \varphi_1$ d'applications $\varphi_1$ et $\varphi_2$, bijectives en dehors de $N_1$ et $N_2$ respectivement, est bijective en restriction au bor\'elien $\varphi_1^{-1}(\varphi_1(X\backslash N_1)\cap X\backslash N_2)$, dont le compl\'ementaire est n\'egligeable). On note $\Def(\Q,\N)$ le quotient de ce groupe obtenu en identifiant deux applications co\"\i ncidant presque s\^urement.

\bigskip

On dit que $(\Q,\N)$ est \textit{ergodique} si toute partie bor\'elienne de $\Q$ est n\'egligeable ou de compl\'ementaire n\'egligeable. Notons que dans le cas ergodique $\B(T)$ co\"\i ncide avec $T$ aux parties n\'egligeables pr\`es. Suivant le point de vue \'evoqu\'e en introduction, une propri\'et\'e d'un espace singulier ergodique $\Q$ est une propri\'et\'e de la notion de quasi-p\'eriodicit\'e choisie. On dira ainsi  qu'une telle propri\'et\'e est un \textit{param\`etre de quasi-p\'eriodicit\'e} (et qu'un invariant d'isomorphisme est une \textit{constante de quasi-p\'eriodicit\'e}).

\bigskip 

Soit $(\Q,\N)$ un espace mesur\'e singulier ergodique. Il existe sur $\Q$ au plus une (\`a constante multiplicative pr\`es) mesure transverse invariante $\sigma$-finie $\Lambda$ d\'efinissant $\N$. De plus, suivant les valeurs que peut prendre $\Lambda$, $(\Q,\N)$ peut \^etre de l'un des trois types suivants :

\medskip

\begin{tabular}{l l l}
\quad - type $\I$  &: & $\Im \Lambda = \{0,\lambda,2\lambda,\ldots,\infty\}=\lambda\NI\cup \{\infty \}$, o\`u $\lambda>0$,\\
\quad - type $\II$ & : & $\Im \Lambda = [0,\infty]$,\\
\quad - type $\III$ & : & $\Im \Lambda = \{0,\infty\},$\\
\end{tabular}

\medskip

\noindent o\`u, pour les espaces de type $\III$, toute mesure transverse invariante $\Lambda$ est triviale (ne contient aucune information autre que $\N$). Il existe \`a isomorphisme pr\`es une unique d\'e\-sin\-gu\-la\-ri\-sa\-tion discr\`ete proprement infinie de $\Q$ ; plus pr\'ecis\'ement, les d\'e\-sin\-gu\-la\-ri\-sa\-tions discr\`etes d'un espace de type $\III$ sont toutes conjugu\'ees, et les d\'e\-sin\-gu\-la\-ri\-sa\-tions discr\`etes d'un espace de type $\II$  sont classifi\'ees (\`a conjugaison pr\`es) par leur mesure transverse. Les d\'efinitions pr\'ec\'edentes s'\'etendent de fa\c con naturelle au cas non ergodique, et tout espace mesur\'e singulier admet une d\'ecomposition $\Q=\Q_\I \amalg \Q_\II \amalg \Q_\III$ en composantes de chaque type, unique aux parties n\'egligeables pr\`es. Ces r\'esultats sont de Murray et von Neumann.

\bigskip
\bigskip

\textit{Convention.} Au cours de ce texte, il n'est question que d'espaces singuliers munis d'une famille $\N$ de bor\'eliens n\'egligeables (i.e. d'espaces mesur\'es singuliers), et on omettra d\'esormais de pr\'eciser cet ensemble dans les notations. On omettra \'egalement l'adjectif \og mesur\'e\fg\ pour qualifier les espaces singuliers.

\bigskip

Consid\'erons un espace singulier ergodique $\Q= (\Q,\N)$ et notons $\Def(\Q) = \Def(\Q,\N)$ son groupe d'automorphismes. Supposons que $\Q$ soit de type $\II$ et fixons une mesure transverse invariante $\sigma$-finie $\Lambda$.

\medskip

\begin{lem}  Consid\'erons $\rho \in \Def(\Q)$. Il existe un unique nombre $\lambda \in ]0,\infty[$ tel que, pour tout isomorphisme d\'esingularisant $\overline \rho : X \to X'$ entre deux d\'e\-sin\-gu\-la\-ri\-sa\-tions discr\`etes, on a $\Lambda(\overline \rho(\Omega))=\lambda\Lambda(\Omega)$ pour tout bor\'elien $\Omega \subset X$. 
\end{lem}

\begin{dem}
Quitte \`a remplacer $\Q$ par $ \Q\backslash N$, o\`u $N\subset \Q$ est n\'egligeable, on peut supposer que $\rho$ est bijective. Consid\'erons deux d\'e\-sin\-gu\-la\-ri\-sa\-tions bijectives $\rho_1 : X_1 \to X_1'$ et $\rho_2 : X_2 \to X_2'$ de $\rho$, o\`u $p_i : X_i \to \Q$ et $p_i' : X_i' \to \Q$ sont des d\'esingularisations discr\`etes de $\Q$ ($i=1,2$). La mesure transverse invariante $\Lambda$ \'etant unique \`a un facteur multiplicatif pr\`es, il existe deux nombres $\lambda_1$ et $\lambda_2$ tels que  $\Lambda(\rho_i(\Omega))=\lambda_i\Lambda(\Omega)$ pour tout bor\'elien $\Omega \subset X_i$. Par d\'efinition il existe deux applications $\varphi : X_1 \to X_2$ et, respectivement, $\varphi' : X_1'\to X_2'$ telles que $p_2\varphi=p_1$ et $p_2'\varphi'=p_1'$. D'apr\`es le lemme \ref{SOE}, $\varphi$ et, respectivement, $\varphi'$, sont en restriction \`a des bor\'eliens non n\'egligeables $A\subset X_1$ et $A'\subset X'_1$, des isomorphismes entre les relations restreintes ${\R_{p_1}}_{|A}$ et ${\R_{p_2}}_{|B}$, et, respectivement, ${\R_{p_1'}}_{|A'}$ et ${\R_{p_2'}}_{|B'}$, o\`u $B=\varphi(A)$ et $B'=\varphi(A')$. Notons que  $\varphi$ pr\'eserve la mesure $\Lambda$ (invariance de $\Lambda$ par conjugaison). De plus, en conjuguant chacune des \'equivalences stables par des isomorphismes partiels int\'erieurs, on peut supposer que $A = X_1$ ou $B=X_2$ (resp. $A' = X_1'$ ou $B'= X_2'$). Supposons par exemple $A=X_1$ et $B'=X_2$. On a alors $\rho_1 \circ \varphi^{-1} = {\varphi'}^{-1} \circ \rho_2$ sur le bor\'elien non n\'egligeable $\varphi(X_1)\subset X_2$. Donc $\lambda_1=\lambda_2$.
\end{dem}

\bigskip

L'application $\mod_\Lambda : \Def(\Q) \to \RI_+^*$ qui \`a $\rho$ associe $\lambda$ est un morphisme de groupes. On note $\Def_\Lambda(\Q)$ son noyau et $F_\Lambda(\Q) \subset \RI_+^*$ son image (groupe fondamental de $\Q$). On a donc une suite exacte 
\[
\begin{CD}
1 @> >> \Def_\Lambda(\Q) @> >> \Def(\Q) @> \mod_\Lambda{} >> \F_\Lambda(\Q) @> >> 1.
\end{CD}
\]
Notons que $\F(\Q)=\F_\Lambda(\Q)$ ne d\'epend \`a multiplication par un scalaire strictement positif  pr\`es que de la classe d'isomorphisme de $\Q$.

\bigskip

Le groupe fondamental d'un espace singulier hyperfini est $\RI_+^*$. Damien Gaboriau a donn\'e de nombreux exemples d'espaces singuliers \`a groupe fondamental trivial en faisant usage des nombres de Betti $L^2$ pour les relations d'\'equivalence \cite{Gaboriau02} (ou alternativement du co\^ut). Le $r$-i\`eme nombre de Betti $L^2$ de $(\Q,\Lambda)$  est le nombre r\'eel positif 
\[
\beta_r(\Q,\Lambda) = \Lambda(X) \cdot \beta_r(\R_p,\Lambda_1),
\]
o\`u $p : X \to \Q$ est une d\'e\-sin\-gu\-la\-ri\-sa\-tion discr\`ete de type $\IIi$ (i.e. telle que $\Lambda(X)<\infty$), $\Lambda_1 = \Lambda/\Lambda(X)$ est la mesure de probabilit\'e sur $X$ associ\'ee \`a $\Lambda$, et $\beta_r(\R_p,\Lambda_1)$ est le $r$-i\`eme nombre de Betti de la relation d'\'equivalence $\R_p$, d\'efini dans \cite{Gaboriau02} (observons que $\beta_r(\Q,\Lambda)$ ne d\'epend pas de la d\'esingularisation discr\`ete choisie, cf. \cite[Corollaire 5.5]{Gaboriau02}). Les nombres $\beta_r(\R_p)=\beta_r(\R_p,\mu)$, o\`u $\mu=\Lambda_1$ est l'unique mesure de probabilit\'e invariante par $\R_p$, sont invariants par isomorphisme de relation d'\'equivalence mesur\'ee.  Par suite, si les relations d'\'equivalence mesur\'ees $\R_p$ et $\R_p'$ associ\'ees \`a deux d\'e\-sin\-gu\-la\-ri\-sa\-tions discr\`etes $p : X\to \Q$ et $p' : X' \to \Q$ de $\Q$ sont isomorphes, et que $\beta_r(\Q,\Lambda) \neq 0$ pour un indice $r$, alors $\Lambda(X) = \Lambda(X')$, et $F(\Q)$ est trivial. (La suite des nombres de Betti $L^2$ \`a multiplication par un scalaire strictement positif  pr\`es est une constante de quasi-p\'eriodicit\'e au sens ci-dessus.)

\bigskip

\begin{prop} Si $p : X \to \Q$ est une d\'e\-sin\-gu\-la\-ri\-sa\-tion de $\Q$, on a alors 
\[
\Def_\Lambda(\Q) \simeq \Out(\R_p) = \Aut(\R_p)/\Int(\R_p).
\] 
\end{prop}

\begin{dem} Il est facile de voir que l'application $p$ induit un morphisme $\Aut(\R_p)\to \Def_\Lambda(\Q)$  de noyau $\Int(\R_p)$. Ce morphisme est surjectif. En effet soit $\rho : \Q \to \Q$ une bijection pr\'eservant $\Lambda$. Il existe une bijection d\'esingularisante $\overline \rho : \Omega \to \Omega$, o\`u $\Omega$ est une partie bor\'elienne de $X$. On peut alors \'etendre $\overline \rho$ \`a $X$ en consid\'erant des expressions de la forme $\psi_i\overline \rho\phi_i$ o\`u $\phi_i : \Omega_i \to \Omega$ et $\psi_i : \Omega \to \Omega_i$  sont des isomorphismes partiels (int\'erieurs) de $\R_p$ et $X = \Omega \amalg \Omega_i$ est une partition bien choisie.
\end{dem}

\bigskip

Terminons par une d\'efinition.

\begin{dfn}
On dira qu'un espace singulier $\Q$ est \textnormal{de type fini} si toute d\'e\-sin\-gu\-la\-ri\-sa\-tion discr\`ete est de type fini, i.e. poss\`ede un graphage u.l.f. (uniform\'ement localement fini). 
\end{dfn}

Notons qu'un espace singulier est de type fini si et seulement si toute d\'e\-sin\-gu\-la\-ri\-sa\-tion discr\`ete peut \^etre engendr\'ee par un nombre fini d'isomorphismes partiels.

\bigskip
\bigskip

\textit{R\'ef\'erences.} La notion d'\'equivalence orbitale stable a \'et\'e introduite par Mackey \cite{Mackey66}. Elle a depuis \'et\'e \'etudi\'ee des deux points de vue bor\'elien et mesur\'e ; cf. \cite{JKL} et \cite{Furman99,Gaboriau02} pour des r\'ef\'erences r\'ecentes. La consid\'eration d'espaces quotients singuliers et de d\'esingularisations a \'et\'e initi\'ee par Connes dans \cite{Connes79}. Nous renvoyons \'egalement \`a \cite{Connes73}, par exemple, pour d'autres d\'eveloppements.

\vspace{1cm}
\section{Structures quasi-p\'eriodiques et repr\'esentations}\label{QP}

Notre but dans ce paragraphe est de pr\'esenter les notions de repr\'esentation hilbertienne et de structure simpliciale associ\'ees \`a une relation d'\'equivalence mesur\'ee, suivant \cite{Connes79} et \cite{Gaboriau02}. Comme annonc\'e en introduction, nous en profitons pour analyser le concept de quasi-p\'eriodicit\'e, et notamment sa formalisation \`a l'aide de la th\'eorie de la mesure. Une structure quasi-p\'eriodique sur un espace est dans ce formalisme une repr\'esentation de relation d'\'equivalence satisfaisant \`a des propri\'et\'es de type \og fid\'elit\'e\fg. La possibilit\'e nous est offerte de commencer par une description g\'en\'erale de la situation, en termes de cat\'egories et foncteurs, isolant celles des repr\'esentations qui sont susceptibles de conduire aux structures quasi-p\'eriodiques. Cette description est uniquement destin\'ee \`a fixer les id\'ees et reste informelle. Une fois ce contexte g\'en\'eral pr\'ecis\'e, nous pourrons nous pencher plus attentivement sur les cas particuliers des complexes simpliciaux et des espaces de Hilbert. 

\bigskip
\bigskip

\textit{I - Le cadre g\'en\'eral.} Nous suivrons bien entendu le principe de base, de ne consid\'erer que des structures que l'on peut \og contruire effectivement\fg\ (i.e. sans recourir \`a l'axiome du choix). Nous les appelerons \og d\'efinis\-sables\fg\ dans ce paragraphe. La signification exacte de ce terme ne sera pr\'ecis\'ee qu'au paragraphe suivant, concernant les cas particuliers. On le remplacera alors par le terme \og mesurables\fg\ pour \'eviter toute ambigu\"\i t\'e.

\bigskip

Fixons une cat\'egorie $\cal C$ dont les objets sont des ensembles.

\bigskip
Soit $\R$ une relation d'\'equivalence (ou un groupo\"\i de) bor\'elienne sur un espace $X$, consid\'er\'ee comme une petite cat\'egorie dont les objets sont les points de $X$ et les morphismes les \'el\'ements de $\R$. On appelle \textit{repr\'esentation de $\R$ dans $\cal C$} la donn\'ee d'un foncteur \textit{d\'efinissable} $F : \R \to \cal C$ (que l'on supposera contravariant).

\bigskip

Lorsque $\R$ est une relation d'\'equivalence mesur\'ee, on appelera encore repr\'esentation de $\R$ la donn\'ee d'un foncteur $F$ dont les lois de composition ont lieu \textit{presque s\^urement} (i.e. $F$ est une repr\'esentation d'une relation $\R_{|X'}$ o\`u $X'\subset X$ est un bor\'elien de compl\'ementaire n\'egligeable).

\bigskip

Soit $\Q$ un espace singulier et $p:X\to \Q$ une d\'esingularisation discr\`ete. On note $\R=\R_p$ la relation d'\'equivalence mesur\'ee associ\'ee. Une repr\'esentation $F$ de $\R$ induit une action (encore not\'ee $F$)  de $\R$ sur l'ensemble
\[
F(X) = \amalg_{x\in X} F(x).
\]
On consid\`ere pour chaque $q\in\Q$ l'ensemble quotient 
\[
\underline F(q) = \amalg_{x\in p^{-1}(q)} F(x) / \sim,
\]
obtenu en identifiant $F(y)$ et $F(x)$ par $F(x,y)$, comme un objet de $\cal C$. On obtient ainsi un diagramme commutatif
\begin{diagram}[nohug,notextflow]
     X   &  \rTo^{F}     & O({\cal C})  \\
   \dTo    &    &  \dTo       \\
  \Q    & \rTo^{\underline F} & O({\cal C}) 
\end{diagram}
o\`u $\underline F$ est une application d\'efinissable de $\Q$ dans les objets $O({\cal C})$ de $\cal C$.

\bigskip

\textit{Un exemple.} Prenons pour $\cal C$ la droite r\'eelle $\RI$ (sans morphisme), pour laquelle la notion de foncteur d\'efinissable (i.e. mesurable) $F :\R \to \RI$ est claire. Tout foncteur mesurable est, par ergodicit\'e, presque s\^urement constant, et l'espace $\underline F(\Q)$  est, \`a un n\'egligeable pr\`es, un nombre r\'eel usuel. Cet exemple concerne plus g\'en\'eralement toute cat\'egorie dont les objets sont des points, i.e. lorsqu'on choisit pour $\cal C$ un espace topologique s\'eparable sans morphisme : un point ne poss\`ede pas de notion int\'eressante de quasi-p\'eriodicit\'e (au sens suivant). 

\bigskip

On appelle \textit{\'el\'ement $\Q$-p\'eriodique de $\cal C$} le $\R$-espace constitu\'e par l'image 
\[
F(X)
\]
d'un foncteur $F$ de $\R$ dans $\cal C$, \textit{lorsque que l'action de $\R$ sur $F(X)$ admet un domaine fondamental d\'efinissable} au sens o\`u il existe une partie d\'efinissable
\[
D=\amalg_{x\in X} D_x \subset F(X)
\]
qui rencontre chaque orbite exactement une fois. En d'autres termes on a
\[
F(x) = \amalg_{y\sim x} F(x,y)D_y.
\]
pour presque tout $x\in X$.

\bigskip

\`A un \'el\'ement $\Q$-p\'eriodique de $\cal C$ est associ\'ee la \textit{lamination} $\underline F(\Q)=F(X)/\R_p$ obtenue en consid\'erant l'espace quotient  de $F(X)$  par l'action de $\R_p$ via $F$.

\bigskip 

\textit{Commentaires.} 1. - Ainsi, bien qu'un espace singulier soit un ensemble au sens usuel du terme, ses points contiennent \textit{a priori} trop d'informations pour \^etre v\'eritablement consid\'er\'es comme des points. La structure singuli\`ere sur cet ensemble d\'etermine des relations entre chacun de ses points, et ces relations prennent ensuite effet lorsque l'on substitue \`a chacun d'eux un \'el\'ement d'une cat\'egorie fix\'ee par la proc\'edure ci-dessus. 

2. - On comprend facilement le terme \og quasi-p\'eriodique\fg\ lorsque le domaine fondamental $D$ est \textit{localement trivial}, au sens o\`u il admet une partition d\'enombrable $D=\amalg_i(X_i\times D_i)$, o\`u $X_i\subset X$ est une partie d\'efinissable non n\'egligeable : par ergodicit\'e presque toute feuille de cette lamination contient pour tout $i$ une infinit\'e de copie de $D_i$ \og uniform\'ement r\'eparties\fg\ (chaque partie $D_i$ apparaissant avec une certaine \og proportion\fg\ dans le cadre mesur\'e). Dans le cas g\'en\'eral, il y a une \og d\'ependance d\'efinissable\fg\ entre deux parties $D_x$ et $D_y$ pour $x$ et $y$ dans une m\^eme feuille. Le cas non ergodique concerne le m\'elange de diff\'erents concepts de quasi-p\'eriodicit\'e.

3. - La d\'efinition de domaine fondamental donn\'ee ici n\'ecessitera d'\^etre adapt\'ee \`a chaque cat\'egorie particuli\`ere. Par exemple, pour la cat\'egorie des complexes simpliciaux, on supposera que chaque partie $D_x$ est simpliciale.

4. - La lamination associ\'ee \`a un \'el\'ement $\Q$-p\'eriodique de $\cal C$ s'identifie de fa\c con d\'efinisssable au domaine fondamental $D$. Le choix d'une section d\'efinissable du fibr\'e $D\to X$ d\'etermine alors un plongement transverse de $X$ de cette lamination.

5. -  On reprend essentiellement ici des techniques introduites par A. Connes dans \cite{Connes79}. Rappelons qu'il est construit dans \cite{Connes79}, \`a l'aide de ces techniques, une th\'eorie de l'int\'egration (transverse) \og en pr\'esence d'un groupo\"\i de mesur\'e\fg. Celle-ci permet d'int\'egrer des \og fonctions positives\fg\ d\'efinies sur l'espace $\Q$ des orbites de ce groupo\"\i de, o\`u une fonction positive est une application qui \`a une orbite $q\in Q$ associe un espace mesur\'e standard $(Y_q,\alpha_q)$, o\`u $\alpha_q$ est une mesure positive sur $Y_q$. D\'ecrire la mesurabilit\'e d'une telle application conduit naturellement \`a la notion de foncteur mesurable du groupo\"\i de vers les espaces mesur\'es. L'int\'egration d'une telle fonction s'effectue \`a l'aide d'une mesure transverse quasi-invariante, associ\'ee \`a un cocycle $\delta$ d\'efini sur le groupo\"\i de (son module).

\bigskip
\bigskip

\textit{II - Exemples de cat\'egories.} D\'efinissons plus pr\'ecis\'ement la notion de mesurabilit\'e pour les cat\'egories suivantes : 
\begin{itemize}
\item ensembles d\'enombrables/bijections,
\item espaces bor\'eliens standard/isomorphismes mesurables,
\item complexes simpliciaux/isomorphismes simpliciaux,
\item espace de Hilbert/op\'erateurs unitaires.
\end{itemize}

\bigskip

Soit $\Q$ un espace mesur\'e singulier. 
\bigskip

\textit{1 - La cat\'egorie des ensembles d\'enombrables ou des bor\'eliens standard.} Commen\c cons par un exemple simple. Soit $p: X \to \Q$ une d\'e\-sin\-gu\-la\-ri\-sa\-tion discr\`ete de $\Q$. Le foncteur naturel $F : X \to R_p$, qui \`a $x\in X$ associe $\R_p^x=\{(x,y)\}\subset \R_p$ et \`a $(x,y)$ l'application \'evidente $\R_p^y\to\R_p^x$, est  mesurable et d\'efinit l'ensemble d\'enombrable $\Q$-p\'eriodique $\R_p$. La lamination associ\'ee est l'espace $X=\R_p/\R_p$ et on a $\underline F = p^{-1}$.\\

Plus g\'en\'eralement soit $\R$ une relation d'\'equivalence mesur\'ee et $F : \R \to \cal C$ une repr\'esentation de $\R$ dans la cat\'egorie des espaces bor\'eliens standard et des isomorphismes mesurables (ou des ensembles d\'enombrables et des bijections). On d\'efinit la mesurabilit\'e de $F$ de la fa\c con suivante (cf. \cite{Connes79}). Soit $\tilde \Omega$ la r\'eunion disjointe
\[
\tilde \Omega = \amalg_{x\in X} F(x).
\]
On dit que $F$ est \textit{mesurable} si $\tilde \Omega$ poss\`ede une structure d'espace bor\'elien standard compatible avec les restrictions aux fibres d\'efinie par la projection naturelle  $\tilde p : \tilde \Omega \to X$, telle que $\tilde p$ et l'application  $\R * \tilde \Omega \to \tilde \Omega$ d\'efinie par 
\[
(x,y)*a \mapsto F(x,y)a
\]
soient mesurables (o\`u $\R*\tilde \Omega$ est le produit fibr\'e de $\R$ et $\tilde \Omega$, i.e. l'ensemble des couples $((x,y),a)$ tels que $\tilde p(a)=y$, cf. \cite{Gaboriau02} par exemple). La notion de domaine fondamental \textit{bor\'elien} $D\subset \tilde \Omega$ est claire ($D$ est une partie bor\'elienne rencontrant exactement une fois chaque orbite) et il en r\'esulte une notion d'espace bor\'elien standard $\Q$-p\'eriodique (cf. ci-dessus). \`A chaque espace bor\'elien standard $\Q$-p\'eriodique est associ\'ee une lamination sur l'espace $\Omega = \underline F(\Q)$ (standard) muni de la d\'e\-sin\-gu\-la\-ri\-sa\-tion $p :\Omega \to \Q$ obtenue par passage au quotient de l'application $\tilde p : \tilde \Omega \to X$ par les actions naturelles de $\R$.

\bigskip

\textit{2 - La cat\'egorie des complexes simpliciaux.} Soit $p : \Sigma \to \Q$ une d\'e\-sin\-gu\-la\-ri\-sa\-tion de $\Q$. On dit que $\Sigma$ est \textit{une d\'esingularisation simpliciale de $\Q$} si pour tout $x\in \Sigma$ la classe $\Sigma_x \subset \Sigma$ de $x$ est munie d'une structure simpliciale connexe, i.e. d'une partition $\Sigma_x = \amalg_{i\geqslant 0} \Sigma_x^{(i)}$ en simplexes non orient\'es de dimension $i$ (satisfaisant aux conditions usuelles de compatibilit\'e), telle que les parties $\Sigma^{(i)} \subset \Sigma$ constitu\'ees des simplexes de dimension $i$ soient bor\'eliennes.

Soit $\R$ une relation d'\'equivalence mesur\'ee. On dit qu'un foncteur $F$ de $\R$ dans la cat\'egorie des (r\'ealisations g\'eom\'etriques de) complexes simpliciaux (non orient\'es) est \textit{mesurable} s'il est mesurable en tant que foncteur \`a valeurs dans les espaces bor\'eliens standard, et si les parties $\tilde \Omega^{(i)} \subset \tilde \Omega = F(X)$ constitu\'ees des simplexes de dimension $i$ sont bor\'eliennes. On dit que $F$ admet un domaine fondamental bor\'elien si, de m\^eme que ci-dessus, il existe une partie bor\'elien $D\subset \tilde \Omega$ rencontrant exactement une fois chaque orbite ; la notion de complexe simplicial $\Q$-p\'eriodique en r\'esulte. Elle co\"\i ncide (\`a r\'ealisation g\'eom\'etrique pr\`es) avec la notion de $\R$-complexe simplicial d\'efinie par D. Gaboriau dans \cite{Gaboriau02}.

On associe \`a une d\'esingularisation simpliciale de $\Q$ un complexe simplicial $\Q$-p\'eriodique de la fa\c con suivante. Fixons une d\'esingularisation simpliciale $p:\Sigma \to \Q$ de $\Q$. Notons $X=\Sigma^{(0)}$ la lamination des sommets et $p^{(0)} : X\to \Q$ la d\'esingularisation discr\`ete associ\'ee. Soit $\R = \R_{p^{(0)}}=\tilde \Sigma \cap X\times X$, o\`u $\tilde \Sigma = \Sigma_p \subset \Sigma\times \Sigma$ est la relation d'\'equivalence associ\'ee \`a $p$. Il est facile de voir que le foncteur $F$ de $\R$ dans $\tilde \Sigma_X = (X\times \Sigma) \cap \tilde \Sigma\subset \tilde \Sigma$ qui \`a $x \in X$ associe le complexe $\tilde \Sigma_x\subset \tilde \Sigma_X$ et \`a $(x,y)\in\R$ l'isomorphisme naturel $\tilde \Sigma_y \to \tilde \Sigma_x$ est mesurable. De plus la donn\'ee d'une section mesurable de la projection sur la seconde coordonn\'ee $\tilde \Sigma_X \subset X\times \Sigma \to \Sigma$ (\`a fibres d\'enombrables) d\'etermine un domaine fondamental bor\'elien et, choisissant convenablement cette section, on d\'efinit ainsi un domaine fondamental (bor\'elien) simplicial. R\'eciproquement \`a un complexe simplicial $\Q$-p\'eriodique, donn\'e par un foncteur $F$, on associe la d\'esingularisation simpliciale $\Sigma=F(X)/\R \to \Q$ obtenue en consid\'erant $F$ comme un foncteur \`a valeurs dans la cat\'egorie des espaces bor\'eliens standard, o\`u la partition bor\'elienne $\Sigma=\amalg_{i\geqslant 0} \Sigma^{(i)}$ est donn\'ee par la projection naturelle (\`a fibres d\'enombrables) de $F(X)^{(i)}$ dans $\Sigma$. Les deux op\'erations d\'ecrites dans ce paragraphe sont inverse l'une de l'autre, identifiant ainsi la cat\'egorie des d\'esingularisations simpliciales de $\Q$ \`a celle des complexes simpliciaux $\Q$-p\'eriodiques.

\bigskip

\textit{Remarques.} i.  Un complexe simplicial quasi-p\'eriodique est n\'ecessairement localement trivial (au sens d\'efini au paragraphe I) du fait qu'il ne poss\`ede qu'un nombre d\'enombrable de g\'eom\'etries locales possibles.

ii.  Le groupe des automorphismes $\Aut(\tilde \Sigma)$ d'un complexe simplicial quasi-p\'eriodique $\tilde \Sigma$ est form\'e des bijections bor\'eliennes de $\Sigma=\tilde \Sigma/\R$, non singuli\`eres et d\'efinies \`a un n\'egligeable pr\`es, qui respectent la structure simpliciale longitudinale. Le sous-groupe distingu\'e $\Int(\tilde \Sigma)$ est form\'e des \'el\'ements de $\Aut(\tilde \Sigma)$ qui fixent l'espace singulier quotient. Notons que $\Out(\tilde \Sigma)=\Aut(\tilde \Sigma)/\Int(\tilde \Sigma)$ diff\'ere en g\'en\'eral de $\Out(\R)$.

iii.  \'Etant donn\'e un complexe simplicial $\Q$-p\'eriodique $\tilde \Sigma$, la famille des bor\'eliens standard $\Q$-p\'eriodiques inclus dans $\tilde \Sigma$ forme une sous-tribu de la tribu bor\'elienne de $\tilde \Sigma$. La donn\'ee d'une mesure quasi-invariante $\mu$ sur $X=\Sigma^{(0)}$ et d'un syst\`eme de Haar $s$ sur $\tilde \Sigma$ (i.e. un champ invariant $(s^x)_{x\in X}$ de mesures sur $\tilde \Sigma$ telle que la mesure $s^x$ soit port\'ee par le complexe simplicial $\tilde \Sigma_x\subset \tilde \Sigma$) d\'etermine une mesure $\h_\mu^{(s)}$ sur cette tribu (cf. \cite{Connes79,Gaboriau02}).

\bigskip

\textit{Complexe quasi-p\'eriodique universel.} Il existe un unique (\`a isomorphisme pr\`es) complexe simplicial $\Q$-p\'eriodique, contenant une copie isom\'etrique de tout autre complexe simplicial $\Q$-p\'eriodique \cite{Gaboriau02}.

\bigskip

On dit qu'un complexe simplicial $\Q$-p\'eriodique est un \textit{arbre $\Q$-p\'eriodique} si presque toute ses classes sont des arbres, et qu'un espace singulier $\Q$ est \textit{arborable} s'il existe un arbre $\Q$-p\'eriodique. De m\^eme on d\'efinit ainsi les notions de dimension, $p$-connexit\'e, etc.,  d'un complexe simplicial quasi-p\'eriodique (cf. \cite{Gaboriau99,Gaboriau02}).

\bigskip

\textit{3 - La cat\'egorie des espaces de Hilbert.} Rappelons enfin la notion de foncteur mesurable \`a valeurs dans la cat\'egorie hilbertienne (\cite{Connes79}).

\bigskip

Soit $X$ un espace bor\'elien standard et $\R$ une relation d'\'equivalence bor\'elienne \`a classes d\'enombrables sur $X$. Soit $H$ un champ mesurable d'espaces hilbertiens de base $X$
(cf. \cite{Dixmier69}). 

\bigskip

Une \textit{repr\'esentation unitaire de $\R$ sur $H$} est la donn\'ee d'une famille d'op\'erateurs unitaires
\[
\pi(x,y) : H_y \to H_x,
\]
$(x,y)\in \R$, satisfaisant aux conditions de composition et de mesurabilit\'e suivantes :
\begin{itemize}
\item $\pi(x,x)=\id\ \mathrm{et} \  \pi(x,z)= \pi(x,y)\pi(y,z)$ pour tout $x\sim y \sim z$.
\item les coefficients
\[
(x,y) \mapsto \langle\pi(x,y) \xi_y |  \eta_x \rangle_x
\]
sont mesurables pour tous champs de vecteurs mesurables $\xi, \eta : X \to H$.
\end{itemize}

\bigskip

\textit{Exemples.} \textit{Repr\'esentation triviale de $\R$.}  La repr\'esentation triviale de $\R$  est la famille $(x,y)\mapsto 1 \in \SI^1 \subset \CI$ op\'erant sur le champ constant $H=X\times \CI$ de fibre $\CI$ (cf. \cite{Pichot04_II}). 

\textit{Repr\'esentation r\'eguli\`ere de $\R$.} On consid\`ere le champ d'espaces de Hilbert $H : x\mapsto \ell^2(\R^x,\h^x)$ qui \`a tout point $x\in X$ associe l'espace des fonctions de carr\'e int\'egrable sur la  classe d'\'equivalence de $x$ pour la mesure de d\'ecompte horizontal $\h^x$. Comme $\h^x=\h^y$ pour tout $(x,y)\in\R$, les espaces $\ell^2(\R^y,\h^y)$ et $\ell^2(\R^x,\h^x)$ sont naturellement identifi\'es par un op\'erateur unitaire $\pi(x,y)$ ; explicitement,
\[
\pi(x,y) : \ell^2(\R^y) \to \ell^2(\R^x)
\] 
est d\'efini par $\pi(x,y)f(x,z) = f(y,z)$.\\

\textit{Int\'egration d'une repr\'esentation.} Soit $\pi$ une repr\'esentation unitaire de $\R$ sur un champ d'espaces de Hilbert $H$ de base $X$. Soit $\mu$ une mesure de probabilit\'e quasi-invariante sur $X$. Consid\'erons l'espace de Hilbert $L^2(X,H)$ des sections de $H$ de carr\'e int\'egrable pour $\mu$. On d\'efinit pour tout sous-groupe d\'enombrable $\G \subset [R]$ du groupe plein de $\R$ une repr\'esentation unitaire $\overline\pi$ sur $L^2(X,H)$ par la formule
\[
\overline\pi(\gamma)\xi : x \mapsto \pi(x,\gamma^{-1} x)\xi_{\gamma^{-1} x}\sqrt{\delta(x,\gamma^{-1}x)}
\]
o\`u $\delta$ est le module de $\mu$.

\bigskip

Consid\'erons deux repr\'esentations $\pi$ et $\pi'$ d'une relation d'\'equivalence mesur\'ee $\R$. On appelle \textit{op\'erateur d'entrelacement} entre $\pi$ et $\pi'$ la donn\'ee d'un champ mesurable essentiellement born\'e d'op\'erateurs $(T_x)_{x\in X}$ tels que 
\[
T_x\pi(x,y)=\pi'(x,y)T_y
\]
pour tout $(x,y)\in\R$.

\bigskip
\bigskip

\textit{Espaces de Hilbert quasi-p\'eriodiques.} Nous n'utiliserons pas et ne d\'efinirons pas cette notion dans cet article. Il sera naturel de dire qu'une repr\'esentation hilbertienne poss\`ede domaine fondamental si l'action obtenue par restriction au fibr\'e en sph\`ere unit\'e en poss\`ede un, en un sens \`a pr\'eciser. L'exemple fondamental est celui de la repr\'esentation r\'eguli\`ere associ\'ee \`a une d\'esingularisation discr\`ete $X\to \Q$ d'un espace singulier $\Q$. La lamination associ\'ee est donn\'ee par $[x] \to \ell^2([x])$ qui \`a une classe $[x]\in \Q$ associe l'espace des fonctions de carr\'e int\'egrable d\'efinies sur cette classe. Pour des consid\'erations r\'ecentes sur la dynamique des groupes d\'enombrables de transformations unitaires sur la sph\`ere unit\'e d'un espace de Hilbert, nous renvoyons \`a \cite{Pestov00} par exemple.

\bigskip
\bigskip

\textit{III - Quelques commentaires.}

\bigskip

\textit{P\'eriodicit\'e.} L'espace singulier d\'efinissant les structures p\'eriodiques est le point (l'espace singulier ergodique de type I). Les d\'e\-sin\-gu\-la\-ri\-sa\-tions associ\'ees proviennent d'actions (propres) de groupes localement compacts : il est n\'ecessaire ici d'\'etudier  les groupes d\'esingularisants sans se restreindre aux seules relations d'equivalence. Notons qu'il est \'egalement possible de mener une telle \'etude pour les espaces singuliers non triviaux, o\`u l'on \'etudie des groupo\"ides d\'e\-sin\-gu\-lari\-sants \cite{Connes79}. Pour les besoins cet article, cependant, les relations d'\'equivalence suffisent.

\bigskip

\textit{P\'eriodicit\'e et quasi-p\'eriodicit\'e.} Soit $\G$ un groupe d\'enombrable et $Y$ un $\G$-complexe simplicial cocompact (ainsi $Y$ est un complexe simplicial p\'eriodique). On peut munir $Y$ d'une structure quasi-p\'eriodique de la fa\c con suivante. Soit $\alpha$ une action ergodique de $\G$ sur un espace de probabilit\'e $(X,\mu)$ et $\R=\R_\alpha$ la relation d'\'equivalence associ\'ee. L'action diagonale de $\G$ sur $X\times Y$ d\'etermine une lamination
\[
\Sigma = (X\times Y)/\G,
\]
qui est une d\'esingularisation simpliciale de $\Q=X/\R$. On dira que le complexe simplicial $\Q$-p\'eriodique $\tilde \Sigma$ associ\'e \`a $\Sigma$ \textit{d\'efinit une structure $\Q$-p\'eriodique sur $Y$}. (cf. \cite[\textsection3.4]{Gaboriau02})

Faisons une constatation simple, extraite de \cite[\textsection 5.33]{Gromov99_MS}, illustrant cette id\'ee. D\'e\-signons par $\RI$ la droite r\'eelle. Tout recouvrement p\'eriodique, disons $\ZI$-invariant (o\`u $\ZI$ agit par translation), de cette droite par des intervalles de longueur $n$ a multiplicit\'e au moins $n$. Par ailleurs, on peut construire un recouvrement quasi-p\'eriodique (d\'efini par exemple via un feuilletage irrationnel du tore) de $\RI$ par des intervalles de longueur $n$, avec multiplicit\'e au plus 2.

\bigskip

\textit{D'autres cat\'egories.} Nous \'etudierons prochainement, en collaboration avec S. Vassout, la notion de vari\'et\'e riemannienne quasi-p\'eriodique en relation avec la signature $L^2$. Plus g\'en\'eralement, toute cat\'egorie d'espaces m\'etriques s\'eparables (espaces CAT(0),...) est \textit{a priori} naturellement sujette \`a quasi-p\'eriodicit\'e.

Notons aussi que certains espaces fonctionnels peuvent \'egalement jouir de propri\'et\'es de quasi-p\'eriodicit\'e. Par exemple, si $p : X \to \Q$ est une d\'e\-sin\-gu\-la\-ri\-sa\-tion discr\`ete de $\Q$, le fibr\'e mesurable $x \mapsto \ell^\infty([x])$, o\`u $[x]$ d\'esigne la classe de $x\in X$, d\'efinit \og l'espace $\ell^\infty$ quasi-p\'eriodique\fg\ 
\[
[x] \mapsto \ell^\infty([x]).
\]
Cet espace est muni d'une structure d'alg\`ebre (effectuer les op\'erations classe par classe) et la norme donn\'ee par le supremum essentiel des normes $\ell^\infty$ en fait une alg\`ebre de von Neumann (il s'agit bien s\^ur de $L^\infty(X)$). De m\^eme, si $[x] \mapsto \ell^2([x])$ est l'espace de Hilbert quasi-p\'eriodique associ\'e \`a la repr\'esentation r\'eguli\`ere de $\R_p$, l'espace quasi-p\'eriodique
\[
[x] \mapsto B(\ell^2([x]))
\]
des op\'erateurs born\'es sur $\ell^2([x])$ (o\`u $[x]\mapsto q_{[x]}$ est mesurable au sens o\`u les fonctions $x\mapsto \langle q_{[x]}\xi_x\mid \eta_x\rangle$ sont mesurables pour toutes sections mesurables $\xi,\eta$ de $\amalg_{x\in X}\ell^2([x])$) est une alg\`ebre et le supremum essentiel en fait une alg\`ebre de von Neumann, qui contient $L^\infty(X)$, agissant par multiplication (il s'agit de l'alg\`ebre de von Neumann de la relation d'\'equivalence $\R_p$).  Nous renvoyons ici \`a \cite[I.4.$\gamma$]{Connes95}.

\vspace{1cm}
\section{Ergodicit\'e forte}\label{erg}

Soit $(X,\mu)$ un espace de probabilit\'e. Fixons une famille d\'enombrable $\Phi$ d'isomorphismes partiels de $X$ pr\'eservant la classe de $\mu$ et consid\'erons le pseudo-groupe $\G=\langle \Phi\rangle$ engendr\'e par $\Phi$. On dit qu'une suite $(A_n)_{n\geqslant 0}$ de parties bor\'eliennes de $X$ est \textit{asymptotiquement invariante} sous l'action de $\G$ si pour tout $\Phi$-mot $m$ de domaine $D \subset X$, on a
\[
\mu(m(A_n\cap D)\backslash A_n) \to_n 0.
\]

\bigskip

La notion de suite asymptotiquement invariante est classique et s'\'enonce traditionellement en terme d'actions de groupe. Si $\alpha : \G \to \Aut(X,\mu)$ est un groupe agissant en pr\'eservant la classe de $\mu$, on dit qu'une suite $(A_n)$ de parties bor\'eliennes de $X$ est \textit{asymptotiquement invariante} sous l'action de $\G$ si
\[
\mu(\alpha(\gamma)A_n\Delta A_n) \to_n 0
\]
pour tout $\gamma \in \Gamma$. 

\bigskip

Il est facile de voir que les deux d\'efinitions donn\'ees ci-dessus lorsque $\Phi$ est constitu\'e d'isomorphismes de $X$ co\"incident (voir \'egalement le lemme \ref{SOEAI}). 
Une suite asymptotiquement invariante $(A_n)$ est dite \textit{non triviale} s'il existe $\delta >0$ tel que 
\[
\delta \leqslant \mu(A_n) \leqslant 1- \delta.
\]

\bigskip

Rappelons \'egalement qu'une suite $(A_n)$ est asymptotiquement invariante \'etant donn\'ee une action $\alpha$ d'un groupe  $\G$ si et seulement si 
\[
\mu(\varphi A_n\Delta A_n) \to_n 0.
\]
pour tout $\varphi \in [\R_\alpha]$ et l'existence de suites asymptotiquement invariantes non triviales ne d\'epend donc que de la relation d'\'equivalence $\R= \R_\alpha$ associ\'ee \`a l'action de $\G$. On dit alors qu'une relation d'\'equivalence ergodique $\R$ est \textit{fortement ergodique} si toute suite asymptotiquement invariante sous l'action de $[\R]$ est triviale (cf. \cite{ConnesKrieger77,Schmidt80,ConnesWeiss80,Schmidt81,JonesSchmidt87}). On v\'erifie facilement que cette propri\'et\'e ne d\'epend que de la classe de $\mu$.

\bigskip

\begin{lem} \label{SOEAI} Soit $\G=\langle\Phi\rangle$ un pseudo-groupe d'isomorphismes partiels agissant ergodiquement sur $(X,\mu)$. L'existence de suites asymptotiquement invariantes non triviales sous l'action de $\G$ est invariante par \'equivalence orbitale stable, i.e. est une propri\'et\'e de l'espace singulier $\Q = X/\G$ des orbites de $\G$ sur $X$.
\end{lem}

\begin{dem}
Montrons d'abord que l'invariance asymptotique sous l'action de $\G$ \'equivaut \`a l'invariance asymptotique sous l'action du groupe plein $[\Phi]$ de la relation d'\'equivalence $\R$ engendr\'ee par $\Phi$ (cf. \cite{JonesSchmidt87,HjorthKechris03}). Fixons une suite $(A_n)$ asymptotiquement invariante sous l'action de $\G$. Soient $\varphi \in [\Phi]$ et $\varepsilon >0$ fix\'e. Soient $m_1,\ldots,m_k$ une famille finie de $\Phi$-mots et $\Omega_1,\ldots \Omega_k$ une famille finie de bor\'eliens disjoints tel que $m_i(x)=\varphi(x)$ pour tout $x\in \Omega_i$ et $\mu(\varphi A_\varepsilon) \leqslant \varepsilon/2$, o\`u $A_\varepsilon = X\backslash \amalg \Omega_i$. Alors
\[
\varphi A_n \backslash A_n = \amalg_1^k\ (m_i(A_n\cap \Omega_i)\backslash A_n) \amalg \varphi A_\varepsilon \backslash A_n\subset  \cup_1^k\ (m_i(A_n\cap D_i)\backslash A_n) \cup \varphi A_\varepsilon \backslash A_n
\]
o\`u $D_i$ est le domaine de $m_i$. Pour tout $n$ suffisament grand on a par hypoth\`ese $\mu(m_i(A_n\cap D_i)\backslash A_n) \leqslant \varepsilon/2k$. Donc $\mu(\varphi A_n \backslash A_n)\to 0$. Ainsi 
\[
\mu(\varphi A_n \Delta A_n)\to 0,
\]
car $A_n\backslash\varphi A_n = \varphi(\varphi^{-1}A_n\backslash A_n)$, et $(A_n)$ est asymptotiquement invariante pour $[\Phi]$. 

R\'eciproquement soit $m : D \to D'$ un $\Phi$-mot. Soit $\varphi_1,\ldots,\varphi_k$ une famille d'isomorphismes partiels de domaines respectifs $D_1,\ldots,D_k$ disjoints tels que $m(x)=\varphi_i(x)$ pour tout $x\in D_i$, et tels que $\mu(\varphi_i(D_i)\cap D_i )=0$. On peut supposer $\mu(m(D_\varepsilon))\leqslant \varepsilon$, o\`u $\varepsilon >0$ est un nombre r\'eel fix\'e et $D_\varepsilon =D\backslash \amalg D_i$ \cite{FeldmanMoore77}. Alors
\[
m(A_n\cap D)\backslash A_n =\amalg_1^k\ (\varphi_i(A_n\cap D_i)\backslash A_n) \amalg m(D_\varepsilon)\backslash A_n\subset \cup_1^k\ \overline\varphi_i(A_n)\backslash A_n \cup m(D_\varepsilon)\backslash A_n
\]
o\`u  $\overline \varphi_i$ est l'extension de $\varphi_i$ \`a $X$ en un isomorphisme d'ordre 2 coincidant avec l'identit\'e sur $X\backslash(D_i\cup \varphi_i(D_i))$. Ainsi $\mu(m(A_n\cap D)\backslash A_n)\to 0$, donc $(A_n)$ est asymptotiquement invariante pour $\G=\langle\Phi\rangle$.

Montrons l'invariance par \'equivalence orbitale stable. Soit $Y \subset X$  un bor\'elien non trivial rencontrant toutes les orbites de $\R$. Soit $A_n \subset Y$ une suite asymptotiquement invariante non triviale pour  $\R_{|Y}$. Consid\'erons un graphage $\Phi$ de $R_{|Y}$ et une famille $\Psi$ d'isomorphismes partiels $X\backslash Y \to Y$ de $\R$ dont les domaines forment une partition de $X\backslash Y$. Alors $\tilde\Phi=\Phi \cup \Psi$ est un graphage de $\R$ et la suite $(A'_n)$ des satur\'es de $A_n$ par $\Psi$ est $\tilde \Phi$-asymptotiquement invariante et non triviale (il suffit de v\'erifier l'invariance asymptotique sur les g\'en\'erateurs d'un graphage). 

R\'eciproquement, soit $(A_n)$ une suite asymptotiquement invariante non triviale pour $\R$. Notons qu'il suffit de montrer qu'il existe un bor\'elien $Y_0\subset Y$ rencontrant toutes les orbites de $\R$ contenant des suites asymptotiquement invariantes pour $\R_{|Y_0}$. Consid\'erons un bor\'elien $Y_0\subset Y$ rencontrant toutes les orbites de $\R$ pour lequel il existe une partition $X\backslash Y_0 = \amalg_{i\geqslant 1} Y_i$ en domaines $Y_i$ d'isomorphismes partiels $Y_i \simeq Y_0$. Notons que comme $(A_n)$ est asymptotiquement invariante, on a 
\[
\mu(A_n\cap Y_i)\to_n 0 \ssi \mu(A_n\cap Y_j)\to_n 0
\]
pour tous $i,j\geqslant 0$ fix\'es. Ainsi $(A_n \cap Y_0)_n$ d\'efinit une suite non triviale de $\R_{|Y_0}$, au sens o\`u l'on peut trouver une sous-suite $(A_m')$ de $(A_n \cap Y_0)_n$ dont la mesure converge vers un nombre r\'eel non nul distinct de $\mu(Y_0)$. Cette suite est asymptotiquement invariante car si $\varphi$ est un isomorphisme de $\R_{|Y_0}$, il s'\'etend par l'identit\'e en un isomorphisme de $\R$, et 
\[
\varphi(A_m')\Delta A_m' =  \varphi(A_m)\Delta A_m.
\]
\end{dem}

\bigskip

\begin{dfn} Soit $\Q$ un  espace singulier ergodique. On dit que $\Q$ est \textnormal{fortement ergodique} si toute suite asymptotiquement invariante d'une d\'e\-sin\-gu\-la\-ri\-sa\-tion discr\`ete est triviale. 
\end{dfn}

\bigskip

\begin{thm}[Jones-Schmidt \cite{JonesSchmidt87}] Soit $\Q$ un espace singulier ergodique. On a l'alternative suivante :
\begin{itemize}
\item soit $\Q$ poss\`ede un quotient moyennable ergodique non trivial,
\item soit $\Q$ est fortement ergodique.
\end{itemize}
\end{thm}

\bigskip

Rappelons la d\'efinition suivante, qui s'adapte imm\'ediatement aux espaces singuliers.

\medskip

\begin{dfn}
Soit $\R$, resp. $\underline \R$, une relation d'\'equivalence mesur\'ee ergodique sur un espace de probabilit\'e $(X,\mu)$, resp. $(\underline X, \underline \mu)$. On dit que $\underline R$  \textnormal{est un quotient de} $\R$ s'il existe une application  bor\'elienne surjective non singuli\`ere  $p : X \to \underline X$  tel que $p^{(1)}(\R) = \underline \R$, o\`u $p^{(1)}(x,y) = (p(x),p(y))$.
\end{dfn}

\bigskip

Nous renvoyons pour l'expos\'e des r\'esultats connus concernant l'ergodicit\'e forte, et notamment pour une d\'emonstration du th\'eor\`eme de Jones-Schmidt, \`a l'article r\'ecent de G. Hjorth et A. Kechris \cite[App. 1]{HjorthKechris03}.

\vspace{1cm}
\section{Concentration}\label{concentr}

Ce paragraphe est essentiellement inspir\'e de l'observation suivante de Gromov, extraite de \cite{Gromov00_SQ},\\

\textit{``If $X$ is foliated (i.e. partitioned) into the orbits of an \textnormal{amenable} group $G$ acting on $X$, then the resulting $d$ on $X$ is, essentially, never concentrated. But if $G$ has property T, then it is concentrated.''}\\

Soit un espace m\'etrique-mesur\'e $(X, d, \mu)$ au sens de Gromov \cite{Gromov00_SQ}, o\`u l'on permet que $d(x,y)=\infty$. Plus pr\'ecis\'ement $X$ est un espace bor\'elien standard, $\mu$ est une mesure de probabilit\'e sans atome sur $X$, et $d$ est une application bor\'elienne  satisfaisant aux axiomes traditionnels d'une distance, except\'e que ses valeurs parcourent $[0,\infty]$. De plus, on suppose que, si $\R=\R_d\subset X\times X$ d\'esigne la relation d'\'equivalence bor\'elienne des couples $(x,y)$ de points \`a distance finie, la mesure $\mu$ est $\textit{quasi-invariante}$ relativement \`a $\R$. \\

\textit{Exemple.} Consid\'erons un feuilletage lisse sur une vari\'et\'e compacte. La donn\'ee d'un champ mesurable de m\'etriques sur chaque feuille d\'etermine un espace m\'etrique-mesur\'e, o\`u on pose
\[
d(x,y) = d_\ell(x,y)
\]
si $x,y$ sont sur une m\^eme feuille $\ell$ et $d(x,y) = \infty$ sinon. Notons que la probabilit\'e pour la classe de Lebesgue que deux points soient \`a distance finie est nulle. Ici $\R$ est la relation d'\'equivalence sous-jacente au groupo\"\i de d'holonomie (partition en feuilles).\\

\textit{Exemple.} Soit $\R$ une relation d'\'equivalence mesur\'ee. Un graphage $K$ de $\R$ d\'efinit naturellement une m\'etrique $d_K$ sur les orbites (distance simpliciale).\\

\'Etant donn\'es deux parties bor\'eliennes $A,B\subset X$ on note 
\[
d(A,B)= {\inf_{x\in A}}^\mu\inf_{y\in B} d(x,y)
\]
o\`u le symbole \og $\inf^\mu$\fg\ d\'esigne l'infimum essentiel relativement \`a $\mu$ (et le symbole \og $\inf$\fg\ l'infimum usuel).

\bigskip

\begin{dfn}[\cite{Gromov00_SQ}]
On dit que  $(X, d, \mu)$ est \textnormal{concentr\'e} s'il existe une fonction $c : ]0,1]^2 \to \RI_+$, telle que pour tous bor\'eliens non n\'egligeables $A,B \subset X$, on a
\[
d(A,B) \leqslant c(\mu(A),\mu(B)).
\]
\end{dfn}

\medskip

On v\'erifie facilement que cette d\'efinition est \'equivalente la suivante : pour tout $\delta > 0$ il existe une constante $r_\delta$ telle que pour tous bor\'eliens $A, B \subset X$,  on a
\[
\mu(A),\, \mu(B) \geqslant \delta \impl d(A,B) \leqslant r_\delta.
\]
Notons que la fonction $c(\delta,\delta') = \sup \{ d(A,B) \mid \mu(A)=\delta,\ \mu(B)=\delta'\}$ est d\'ecroissante \`a $\delta$ fix\'e. 

\bigskip

\textit{Remarque.} La d\'efinition qualitative donn\'ee ci-dessus est insuffisante pour obtenir, comme dans la situation classique, des r\'esultats sur la concentration des fonctions 1-lipschitziennes (voir aussi \cite[thm. 4]{Pichot04_II}).

\bigskip

Observons que les espaces m\'etriques-mesur\'es concentr\'es sont en particulier ergodiques, au sens o\`u $\R$ est ergodique relativement \`a $\mu$ ou, de fa\c con \'equivalente, au sens suivant.

\medskip

\begin{dfn}
Soit $(X,d,\mu)$ un espace m\'etrique-mesur\'e au sens ci-dessus. On dit que $d$ est \textnormal{ergodique} (relativement \`a $\mu$) si pour toutes parties  bor\'eliennes $A,B\subset X$ non n\'egligeables, on a $d(A,B)<\infty$. 
\end{dfn}

\medskip

\begin{thm}
Soit $(X,\mu)$ un espace bor\'elien standard et $\Phi$ une famille finie d'isomorphismes partiels (pr\'eservant la classe de $\mu$ et) induisant une m\'etrique ergodique $d=d_\Phi$ sur $X$. Alors
$X$ est concentr\'e si et seulement s'il ne contient pas de bor\'eliens asymptotiquement invariants non triviaux sous l'action du pseudo-groupe $\G=\langle \Phi \rangle$.
\end{thm}

\begin{dem}
Soit $(A_n)$ une suite asymptotiquement invariante telle que
\[
\delta \leqslant \mu(A_n) \leqslant 1-\delta.
\]
Supposons en raisonnant par l'absurde que l'espace $X$ soit concentr\'e et consid\'erons un entier $r=r_{\delta/2}$ telle que si $\mu(A),\, \mu(B) \geqslant \delta/2$ alors $d(A,B) \leqslant r$. Soit $F = \{m^1, m^2, \ldots,m^f\}$ la famille finie des $\Phi$-mots de longueur $\leqslant r$. Notons $B_n=X\backslash A_n$ et fixons un mot $m^i\in F$, de domaine $D^i$. Soit $(D_n^i)$ la suite de bor\'eliens d\'efinis par $D_n^i = \{x\in D^i\cap A_n|\, m^i(x) \in B_n\}$. Par d\'efinition de $A_n$ on a 
\[
\mu(m^i(A_n\cap D^i) \backslash A_n) \to 0,
\]
i.e., $\mu(m^i(D_n^i)) \to 0$, et il existe un entier $N_i$ suffisament grand pour que  pour $n\geqslant N_i$,
\[
B_n^{m_i} = B_n \backslash m^i(D_n^i)
\]
soit de mesure $\geqslant\delta - \frac \delta {2f}$. Choisissons $N = \max_i N_i$ et consid\'erons le bor\'elien $C_N = \cap_i B_N^{m^i}$, de mesure $\geqslant\delta/2$. On a 
\[
\mu( m(A_N \cap D)\cap C_N) = 0
\]
pour tout mot $m\in F$ de domaine $D$. En d'autres termes $d(A_N,C_N) > r$, contrairement \`a l'hypoth\`ese.\\

R\'eciproquement supposons que $X$ ne soit pas concentr\'e et construisons une suite asymptotiquement invariante. Soit $(A_n), \, (B_n)$ deux suites de bor\'eliens de taille $\geqslant\delta$ tels que
\[
d(A_n,B_n) \geqslant n+1.
\]
Ainsi la \og boule\fg\ de centre $A_n$ et de rayon $n$ est disjointe de $B_n$ ; on supposera que $B_n$ co\"\i ncide avec le compl\'ementaire de cette boule. Construisons une suite de fonctions 
\[
\pi_n\in L^\infty(X,\mu)
\]
de la fa\c con suivante. Sur $A_n$ (resp. $B_n$), on pose $\pi_n = 1$ (resp. 0). Sur la sph\`ere de centre $A_n$ et de rayon $i=1\ldots n$, on pose $\pi_n = 1-1/i$. 
Soit $m$ un $\Phi$-mot de domaine $D$ et de longueur $l$. Il est clair que pour presque tout $x \in D$ on a
\[
\abs{\pi_n(m(x)) - \pi_n(x)} \leqslant \frac l n.
\]
Plus g\'en\'eralement pour tout automorphisme partiel $\varphi \in [[\R]]$ tel que 
\[
\abs \varphi  = \int_{D'} d(\varphi^{-1} x,x) d\mu(x) < \infty
\]
o\`u $D'$ est l'image de $\varphi$, on a
\begin{eqnarray}
\norm{(\pi_n \circ \varphi^{-1} - \pi_n)_{|D'}}_1 &=& \int_{D'}\abs{\pi_n(\varphi^{-1} x)-\pi_n(x)}d\mu(x)\nonumber\\
&\leqslant&\int_{D'}\frac {d(\varphi^{-1} x,x)} n d\mu(x)=\frac{\abs \varphi} n \to_n 0.\nonumber
\end{eqnarray}
Notons $A_n^\alpha = \{ \pi_n \geqslant \alpha\}$. On a pour tout $x \in D'$,
\[
\abs{\pi_n (\varphi^{-1} x) - \pi_n(x)} = \sum_{i=0}^n \frac 1 n \chi_{\varphi (A_n^{1-i/n}\cap D)\, \Delta\, (A_n^{1-i/n}\cap D')}(x),
\]
o\`u $D$ est le domaine de $\varphi$. Donc
\[
\norm{(\pi_n \circ \varphi^{-1} - \pi_n)_{|D'}}_1 = \int_0^1 \mu(\varphi (A_n^\alpha\cap D) \Delta (A_n^\alpha\cap D'))d\alpha.
\]
Il en r\'esulte, quitte \`a extraire une sous-suite, que
\[
\mu(\varphi (A_n^\alpha\cap D) \backslash A_n^\alpha)\to 0
\]
pour presque tout $\alpha$.

Par d\'efinition $\abs \varphi =\mu(D')$ pour tout isomorphisme partiel $\varphi \in \Phi$ d'image $D'$. Comme $\Phi$ est de cardinal fini, on obtient ainsi par extractions successives une sous-suite $(A_m)$ de $(A_n)$ telle que $\mu(\varphi(A_m^\alpha \cap D)\backslash A_m^\alpha)\to 0$ pour tout $\varphi\in \Phi$ et presque tout $\alpha$. Comme il suffit de v\'erifier l'invariance asymptotique sur un syst\`eme g\'en\'erateur de $\R$, presque toute suite $(A_m^\alpha)$ est asymptotiquement invariante sous l'action de $\G=\langle \Phi\rangle$. De plus, elles sont non triviales pour $\alpha>0$, et le th\'eor\`eme en r\'esulte.
\end{dem}

\bigskip

\textit{Remarque.} Le th\'eor\`eme pr\'ec\'edent est clairement faux lorsque la m\'etrique $d_\Phi$ n'est pas de type fini (on peut alors  choisir pour $\Phi$ une partition de $\R$ en isomorphismes partiels). Dans ce cas, il est facile d'adapter la d\'emonstration ci-dessus pour \'etablir le r\'esultat suivant. \'Etant donn\'ee une famille d\'enombrable $\Phi$ d'isomorphismes partiels de $X$ pr\'eservant la classe de la mesure $\mu$, \textit{l'action sur $X$ du pseudo-groupe engendr\'e par $\Phi$ contient une suite asymptotiquement invariante non triviale si et seulement si pour toute m\'etrique de type fini $d \geqslant d_\Phi$, l'espace m\'etrique-mesur\'e $(X,\mu,d)$ n'est pas concentr\'e.}

\bigskip

 Soit $\R$ une relation d'\'equivalence mesur\'ee de type fini. On dit que $\R$ est \textit{concentr\'ee} si pour tout graphage fini $\Phi$ de $\R$ l'espace m\'etrique-mesur\'e $(X,d_\Phi,\mu)$ est concentr\'e. Notons que si $\R$ n'est pas concentr\'e, alors aucun des espaces $(X,d_\Phi,\mu)$ associ\'e \`a un graphage fini $\Phi$ ne l'est ; la propri\'et\'e de concentration est ainsi ind\'ependante de la m\'etrique (de type fini) choisie sur l'espace ambiant. On dit qu'un espace singulier ergodique de type fini est \textit{concentr\'e} si toute d\'e\-sin\-gu\-la\-ri\-sa\-tion discr\`ete l'est.  

\bigskip

\begin{thm} Soit $\Q$ un espace singulier ergodique de type fini. Alors $\Q$ est fortement ergodique si et seulement s'il est concentr\'e.
\end{thm}

\bigskip

Rappelons que d'apr\`es le r\'esultat de Schmidt-Connes-Weiss (cf. \cite{Schmidt80,ConnesWeiss80}), un groupe d\'enombrable poss\`ede la propri\'et\'e T de Kazhdan si et seulement si toutes ses actions ergodiques de type $\IIi$ sont fortement ergodiques.

\begin{cor}
Soit $\G$ un groupe de type fini. Alors $\G$ poss\`ede la propri\'et\'e T de Kazhdan si et seulement si toutes ses actions ergodiques de type $\IIi$ sont concentr\'ees.
\end{cor}

\bigskip

Nous renvoyons le lecteur \`a \cite{Gromov00_SQ} pour un \'eventail d'id\'ees sur la concentration. On dit classiquement qu'un espace concentr\'e est, en g\'en\'eral, un espace \og de grande dimension\fg. On retrouve dans le contexte quasi-p\'eriodique ce type ph\'enom\`ene. Par exemple, les shifts de Bernouilli (non moyennables), i.e. $\{0,1\}^\G$ o\`u $\G$ agit par translation, sont concentr\'es (cf. \cite{LosertRindler81,JonesSchmidt87,HjorthKechris03}).

\vspace{1cm}
\section{In\'egalit\'es isop\'erim\'etriques}\label{isop}

Soit $\R$ une relation d'\'equivalence mesur\'ee de type fini sur un espace de probabilit\'e $(X,\mu)$. Soit $K$ un graphage u.l.f.  de $\R$.

\begin{dfn} On dit que $K$ poss\`ede des \textnormal{suites de F\o lner \'evanescentes} (vanishing F\o lner sequences)  relativement \`a $\mu$ s'il existe une suite $(A_n)$ de bor\'eliens non n\'egligeables de $X$ et une suite $(\varepsilon_n)$ de nombres r\'eels convergeant vers 0  telles que  
\[
\mu(A_n) \to 0 \qquad \mathrm{et} \qquad \mu(\del_{K} A_n) \leqslant \varepsilon_n \mu(A_n).
\]
\end{dfn}

\bigskip

\textit{Remarque.} L'existence de telles suites de F\o lner est une propri\'et\'e de $K$ et non de $\R$ ; le r\'esultat principal de ce paragraphe caract\'erise les relations d'\'equivalence dont tout graphage contient des suites de F\o lner \'evanescentes, dans l'esprit de \cite[page 443]{CFW81} (cf. \'egalement \cite{Kaima97}). Notons que dans la d\'efinition pr\'ec\'edente les composantes connexes de $A_n$ (relativement \`a la structure simpliciale) n'ont pas de raison d'\^etre finies (c'est l\`a une diff\'erence essentielle avec les suites de F\o lner de \cite{CFW81}). La notion de suites de F\o lner \'evanescentes est une reformulation g\'eom\'etrique de la notion dynamique de \textit{I-suites} consid\'er\'ee dans \cite{JuncoRosenblatt79,Schmidt81,Rosenblatt81}, lorsque $\mu$ est une mesure de probabilit\'e \textit{invariante}.\\

\textit{Exemple.} Soit $\R$ une relation d'\'equivalence de type $\IIi$. Soit $A \subset X$ un bor\'elien de mesure 1/4. Consid\'erons un graphage born\'e $K_A$ de $\R_{|A}$, une partition infinie $A_1^1, A_2^1, \ldots$ de $A$ et une partition
\[
\{A_i^j\}_{i\geqslant 1,\ 2\leqslant j\leqslant n_i}
\]
de $X\backslash A$ de sorte que $\mu(A_i)=\mu(A_i^j)$ et que $n_i \to_i \infty$. En choisissant des isomorphismes partiels $A_i^j \to A_i^{j+1}$, il est facile de compl\'eter $K_A$ en un graphage born\'e $K$ de $\R$. Ce graphage contient des suites de F\o lner \'evanescentes. $X$ peut cependant \^etre fortement ergodique, si par exemple $\R$ poss\`ede la propri\'et\'e T de Kazhdan.

\bigskip

\textit{Exemple.} Schmidt \cite{Schmidt81} a construit un \textit{arbre} quasi-p\'eriodique de type $\II_1$ et de valence 6, \`a la fois concentr\'e et contenant des suites de F\o lner au sens ci-dessus. Nous avons vu dans l'exemple pr\'ec\'edent que tout espace singulier poss\`ede un graphe quasi-p\'eriodique contenant des suites de F\o lner. Hjorth et Kechris \cite{HjorthKechris03} ont montr\'e que pour tout espace singulier obtenu par action m\'elangeante de type $\IIi$ du groupe libre \`a deux g\'en\'erateurs, on peut construire un arbre quasi-p\'eriodique de valence 4 contenant des suites de F\o lner. Rappelons que le groupe libre \`a deux g\'en\'erateurs permet de d\'efinir une infinit\'e non d\'enombrable de notions de quasi-p\'eriodicit\'e distinctes, en vertu d'un th\'eor\`eme r\'ecent de Gaboriau-Popa \cite{GaboriauPopa04}.

\bigskip

\begin{thm}
Soit $\R$ une relation d'\'equivalence ergodique de type fini pr\'eservant une mesure de probabilit\'e $\mu$. Alors $\R$ poss\`ede un quotient moyennable si et seulement si chacun de ses graphages u.l.f. contient des suites de F\o lner \'evanescentes.
\end{thm}

\medskip

\begin{dem}  D'apr\`es le th\'eor\`eme de Jones-Schmidt, nous devons montrer qu'une relation d'\'equivalence $\IIi$ de type fini est fortement ergodique si et seulement si l'un de ses graphages u.l.f. ne contient pas de suites de F\o lner \'evanescentes.

\medskip

Soit $\R$ une relation d'\'equivalence $\IIi$ de type fini. Supposons que tout graphage u.l.f. de $\R$ contient des suites de F\o lner \'evanescentes et construisons une suite asymptotiquement invariante non triviale. Soit $K$ un graphage u.l.f. de $\R$ et $\varepsilon >0$. On consid\`ere l'ensemble $\cal E$ des bor\'eliens $A$  de $X$ tels que
\[
\mu(\del_K A) \leqslant \varepsilon \mu(A)
\]
et
\[
\mu(A) \leqslant c,
\]
o\`u $c\in ]0,1[$ est fix\'e. \'Etant donn\'es $A$ et $B$ dans $\cal E$ on pose
\[
A\leqslant B\ \mathrm{si}\ A\subset B\quad \mathrm{et}\quad \mu(\del_K A'\backslash A) \leqslant \varepsilon \mu(A'),
\]
o\`u $A' = B\backslash A$ et l'inclusion $A\subset B$ a lieu \`a une partie n\'egligeable pr\`es. Ceci d\'efinit un ordre partiel sur $\cal E$. Soit $\cal E'$ un sous-ensemble totalement ordonn\'e et $A_n \in \cal E'$ une suite telle que $\sup_n \mu(A_n) = \sup_{A\in \cal E'} \mu(A)$. Posons $A_\infty = \cup_n A_n$ et $A_i'=A_i\backslash A_{i-1}$. Alors
\[
\mu(\del_K A_\infty) =\lim_n \mu(\del_K A_n\backslash A_\infty) \leqslant \varepsilon \lim_n\mu(A_n) = \varepsilon \mu(A_\infty),
\]
et, si $A\in {\cal E}'$, on a
\[
\mu(\del_K A_\infty'\backslash A) =\lim_n \mu(\del_K A_n'\backslash A_\infty) \leqslant \varepsilon \lim \mu(A'_n)=\varepsilon \mu(A_\infty')
\]
o\`u $A_\infty'=A_\infty\backslash A$ et $A_n'=A_n\backslash A$. Ainsi $\cal E$ est inductif. Soit $A\in \cal E$ un \'el\'ement maximal (lemme de Zorn). Supposons par l'absurde que $\delta = c - \mu(A) >0$ et notons $K'=X\backslash A \times X\backslash A \cap K$. 

Compl\'etons $K'$ en un graphage u.l.f. $K''$ de $\R_{|X\backslash A}$, et choisir une famille finie d'isomorphismes partiels $\varphi_i : A \to X\backslash A$ de $\R$ dont les domaines partitionnent $A$. Alors 
\[
\overline K = K''\cup_i (\graph\varphi_i \cup \graph\varphi_i^{-1})
\]
est un graphage u.l.f. de $\R$. 

Soit $A_n\subset X$ une suite de F\o lner \'evanescentes pour $\overline K$. Posons $A_n' = A_n \cap X\backslash A$. Par d\'efinition de $\overline K$ on a $\mu(A_n') >0$ pour tout $n$ suffisament grand. Soient
\begin{eqnarray}
B_n^1 &= &\{x\in A_n\cap A |\ \varphi_i(x)\notin A_n'\cup \del_{K''}A_n'\},\nonumber\\
B_n^2 &=& \{x\in A_n\cap A |\ \varphi_i(x)\in \del_{K''}A_n'\},\nonumber\\
B_n^3 &=& \{x\in A_n\cap A |\ \varphi_i(x)\in A_n'\}.\nonumber
\end{eqnarray}
On a
\[
\mu(\del_{\overline K} A_n) \geqslant \mu(\del_{\overline K} B_n^1) + \mu(\del_{K''} A_n').
\]
Par ailleurs $\overline K$ \'etant u.l.b. (rappelons que $\mu$ est invariante), il existe une constante $C=C(\overline K)$ telle que
\[
\mu(B_n^1)\leqslant C\mu(\del_{\overline K} B_n^1), 
\quad \mu(B_n^2)\leqslant C\mu(\del_{K''} A_n'), \quad \mathrm{et} \quad\mu(B_n^3) \leqslant C \mu(A_n').
\]
Par suite pour tout $\varepsilon'$ il existe $n$ tel que
\[
\mu(\del_{\overline K} B_n^1)+\mu(\del_{K''}A_n')\leqslant \varepsilon'(\mu(B_n^1)+\mu(B_n^2)+\mu(B_n^3)+\mu(A_n')),
\]
donc
\[
(1-\varepsilon'C)\mu(\del_{K''} A_n') \leqslant(\varepsilon'C-1) \mu(\del_{\overline K} B_n^1) + \varepsilon'(1+C)\mu(A_n'),
\]
et si $\varepsilon'$ v\'erifie
\[
\varepsilon' C < 1 \quad \mathrm{et} \quad \frac{\varepsilon'(1+C)}{1-\varepsilon'C} \leqslant \varepsilon,
\]
on obtient un entier $n$ tel que que $0< \mu(A_n') <\delta$ et 
\[
\mu(\del_{K''} A_n') \leqslant \varepsilon \mu(A_n').
\]
Alors
\[
\mu(\del_K A_n'\backslash A) =\mu(\del_{K'} A_n')\leqslant \varepsilon \mu(A_n')
\]
donc $A < A \amalg A_n'$. Enfin,
\[
\mu(\del_K(A\amalg A_n'))= \mu(\del_K A\backslash A_n')+\mu(\del_K A_n'\backslash A) \leqslant \varepsilon \mu(A)+\varepsilon \mu(A_n')=\varepsilon \mu(A\amalg A_n'),
\]
donc $A \amalg A_n' \in \cal E$, d'o\`u une contradiction. Ainsi on peut trouver une famille $A_n \subset X$ de bor\'elien tels que $\mu(A_n) = c$ et $\mu(\del_K A_n) \leqslant 1/n$ ; il en r\'esulte, en partionnant $K$ en un nombre fini d'isomorphismes partiels, que $\R$ contient des suites asymptotiquement invariantes (de mesure $c$).\\

R\'eciproquement supposons que $\R$ ne soit pas fortement ergodique et montrons que tout graphage u.l.f. $K$ contient des suites de F\o lner \'evanescentes. (Nous nous inspirons ici d'un argument classique, cf. \cite{JuncoRosenblatt79,Schmidt81,Rosenblatt81,HjorthKechris03}.) Soit $K$ un graphage u.l.f. de $\R$. Pour tout nombre r\'eel $0<\delta<1$ il existe un suite asymptotiquement invariante $A_k^\delta$ telle que $\mu(A_k^\delta)=\delta$ (cf. \cite{JonesSchmidt87}). Soit $n$ fix\'e et $\delta=1/n$. $K$ \'etant u.l.f., on peut le partionner en une famille finie $F$ d'isomorphismes partiels. Soit $k=k(n)$ un entier suffisament grand pour que 
\[
\sum_{\varphi\in F}\mu(\varphi(A_k^\delta\cap D_\varphi) \backslash A_k^\delta) \leqslant \frac 1 {n^2},
\]
o\`u $D_\varphi$ est le domaine de $\varphi$. On pose $A_n' = A_{k(n)}^\delta$. Alors
\[
\mu(A_n')\to 0\quad \mathrm{et}\quad \mu(\del_K A_n') \leqslant \frac 1 n \mu(A_n').
\] 
\end{dem}

\bigskip

\textit{Remarque.}  Il existe des relations relativement simples entre concentration et in\'ega\-lit\'es isop\'erim\'etriques dans le cadre standard \cite{Ledoux01}. Le r\'esultat pr\'ec\'edent est de nature diff\'erente, notamment en ce qu'il concerne seulement des in\'egalit\'es isop\'erim\'etriques \textit{\'evanescentes}.

\bigskip

\'Enon\c cons d\`es maintenant le corollaire suivant (cf. \cite{Pichot04_II}).

\begin{cor}
Soit $\Q$ un espace singulier ergodique ayant la propri\'et\'e T de Kazhdan. Il existe un graphe $\Q$-p\'eriodique u.l.f. ne contenant pas de suites de F\o lner \'eva\-nescentes.
\end{cor}

\begin{dem} Ayant la propri\'et\'e T, un tel espace est de type $\II$, et admet, de plus, une d\'e\-sin\-gu\-la\-ri\-sa\-tion discr\`ete $\IIi$ de type fini. Il existe alors un graphage u.l.f. de cette d\'e\-sin\-gu\-la\-ri\-sa\-tion ne contenant pas de suites de F\o lner \'evanescentes, en vertu du th\'eor\`eme ci-dessus.
\end{dem}

\bigskip

\begin{thm}
Soit $\Q$ un espace singulier ergodique de type fini. Alors $\Q$ poss\`ede un quotient moyennable si et seulement si tout graphe $\Q$-p\'eriodique u.l.f. contient des suites de F\o lner \'evanescentes (relativement \`a une mesure de probabilit\'e dont la classe est d\'etermin\'ee par $\Q$).
\end{thm}

\begin{dem} Le fait que tout graphage u.l.f. contient des suites de F\o lner lorsque $\Q$ poss\`ede un quotient moyennable r\'esulte sans changement de la d\'emonstration ci-dessus. Ainsi, il reste seulement \`a montrer que,  lorsque $\Q$ est de type $\III$ et tel que tout graphe $\Q$-p\'eriodique sym\'etrique u.l.f. contienne des suites de F\o lner \'evanescentes,  il existe une d\'esingularisation discr\`ete de $\Q$ admettant des suites asymptotiquement invariantes non triviales. Consid\'erons un graphage u.l.f. $K$ d'une d\'e\-sin\-gu\-la\-ri\-sa\-tion discr\`ete $\R$ de $\Q$. Soit $\varepsilon >0$.  Reprenons la d\'emonstration du th\'eor\`eme pr\'ec\'edent. Notons $\cal E$ l'ensemble des bor\'eliens $A$  de $X$ tels que
\[
\mu(\del_K A) \leqslant \varepsilon \mu(A)
\]
et
\[
\mu(A) \leqslant c
\]
o\`u $c\in ]0,1[$. Avec la m\^eme relation d'ordre $\cal E$ est inductif. Soit $A\in \cal E$ un \'el\'ement maximal. Supposons que $\delta = c - \mu(A) >0$ et notons $K'=X\backslash A \times X\backslash A \cap K$.  On peut compl\'eter $K'$ en un graphage u.l.f. $\overline K$ de $\R_{|X\backslash A}$. Soit $A_n'\subset X$ une suite de F\o lner \'evanescentes pour $\overline K$ (relativement \`a la mesure $\mu/\mu(X\backslash A)$). Soit un entier $n$ tel que que $0< \mu(A_n') <\delta$ et 
\[
\mu(\del_{K'} A_n') \leqslant\mu(\del_{\overline K} A_n') \leqslant \varepsilon \mu(A_n').
\]
Alors
\[
\mu(\del_K A_n'\backslash A) =\mu(\del_{K'} A_n')\leqslant \varepsilon \mu(A_n')
\]
donc $A < A \amalg A_n' \in \cal E$, d'o\`u une contradiction. Ainsi on peut trouver une famille $A_n \subset X$ de bor\'elien tels que $\mu(A_n) = c$ et $\mu(\del_K A_n) \leqslant 1/n$ ; il en r\'esulte, en partitionnant $K$ en un nombre fini d'isomorphismes partiels, que $\R$ contient des suites asymptotiquement invariantes.
\end{dem}

\bigskip

\textit{Variation.} Une d\'emonstration similaire \'etablit, dans l'esprit de \cite{CFW81}, qu'une relation d'\'equivalence ergodique $\R$ poss\`ede un quotient moyennable si et seulement si chacun de ses sous-graphes sym\'etriques u.l.b. (i.e. dont les d\'eriv\'ees de Radon-Nikodym sont born\'ees mais qui n'engendrent pas n\'ec\'essairement $\R$) poss\`ede des suites de F\o lner \'evanescentes.

\bigskip
\bigskip
\bigskip
\noindent  Unit\'e de Math\'ematiques Pures et Appliqu\'ees,\\
Unit\'e Mixte de Recherche CNRS 5669\\
46, all\'ee d'Italie,\\
69364 Lyon Cedex 07, France.\\
mpichot@umpa.ens-lyon.fr, pichot@ms.u-tokyo.ac.jp.


\bigskip


\begin{thebibliography}{00}

\bibitem{BarrePichot04_I} Barr\'e S., Pichot M., \og Trivialit\'e du groupe d'automorphismes d'un immeuble triangulaire g\'en\'erique\fg, en pr\'eparation.

\bibitem{BarrePichot04_II} Barr\'e S., Pichot M., en pr\'eparation.

\bibitem{Bellissard} Bellissard J., \og The noncommutative geometry of aperiodic solids\fg, Geometric and Topological Methods for Quantum Field Theory, Villa de Leyva, 2001, World Sci. Publishing, River Edge, NJ, 86-156, \textit{2003}.


\bibitem{CantwellConlon98} Cantwell J., Conlon L., \og Generic leaves\fg, Comment. Math. Helv., 73, 306-336, \textit{1998}.

\bibitem{Champetier00} Champetier C., \og L'espace des groupes de type fini\fg,  Topology  39,  no. 4, 657--680, \textit{2000}.

\bibitem{Connes73} Connes A., \og Une classification des facteurs de type $\III$\fg, Ann. Sci. \'Ecole Normale Sup. 4\`eme S\'erie, 6 fasc 2, 133-252, \textit{1973}.


\bibitem{Connes79} Connes A., \og Sur la th\'eorie non commutative de l'int\'egration\fg, Alg\`ebres d'op\'erateurs, Lecture Notes in Math., 725, 19-143, \textit{1979}.

\bibitem{Connes95} Connes A., \og Non commutative geometry\fg, Academic Press, Inc., San Diego, CA, \textit{1994}.

\bibitem{Connes04} Connes A., \og On the foundation of noncommutative geometry\fg, Notes de la conf\'erence NCGOA, Nashville, \textit{2004}.

\bibitem{CFW81} Connes A., Feldman J., Weiss B., \og An amenable equivalence relation is generated by a single transformation\fg, Ergod. Th. \& Dynam. Sys., 1, 431-450, \textit{1981}.

\bibitem{ConnesKrieger77} Connes A., Krieger W., \og Measure space automorphisms, the normalizers of their full groups, and approximate finiteness\fg, J. Functional Analysis, 24, 336-352, \textit{1977}.

\bibitem{ConnesWeiss80} Connes A., Weiss B., \og Property ${\rm T}$ and asymptotically invariant sequences\fg,  Israel J. Math.  37, no. 3, 209--210, \textit{1980}.

\bibitem{Dixmier69} Dixmier J., \og Les alg\`ebres d'op\'erateurs dans l'espace hilbertien (alg\`ebres de von Neumann)\fg, Gauthier-Villars \'Editeurs, Paris, \textit{1969}. Deuxi\`eme \'edition, revue et augment\'ee, Cahiers Scientifiques, Fasc. XXV.


\bibitem{FeldmanMoore77} Feldman J., Moore C.,  \og Ergodic equivalence relations, cohomology, and von Neumann algebras. I,II\fg,  Trans. Amer. Math. Soc., 234(2), 289-359, \textit{1977}.

\bibitem{FHM78} Feldman J., Hahn P., Moore C.C., \og Orbit structure and countable sections for actions of continuous groups\fg, Adv. in Math., 28, 186-230, \textit{1978}.

\bibitem{Furman99} Furman A., \og Gromov's measure equivalence and rigidity of higher rank lattices\fg,  Ann. of Math. (2)  150,  no. 3, 1059--1081, \textit{1999}. 

\bibitem{Gaboriau99} Gaboriau D., \og Co\^ut des relations d'\'equivalence et des groupes\fg,  Invent. Math.  139,  no. 1, 41--98, \textit{2000}.

\bibitem{Gaboriau02} Gaboriau D., \og Invariants $l\sp 2$ de relations d'\'equivalence et de groupes\fg, Publ. Math. Inst. Hautes \'Etudes Sci.  No. 95, 93--150, \textit{2002}.

\bibitem{GaboriauPopa04} Gaboriau D., Popa S., \og An uncountable family of non orbit equivalent actions of $F_n$\fg, pr\'epublication, \textit{2004}.

\bibitem{Ghys95} Ghys \'E., \og Topologie des feuilles g\'en\'eriques\fg, Ann. of Math., 141, 387-422, \textit{1995}.

\bibitem{Gromov99_MS} Gromov M., \og Metric structures for Riemannian and non-Riemannian spaces\fg, Birkha\"user Boston Inc., \textit{1999}.

\bibitem{Gromov00_SQ} Gromov M., \og Spaces and questions\fg, Geom. Funct. Anal., Special Volume (Tel Aviv, 1999), Part I, 118-161, \textit{2000}.

\bibitem{Heafliger58} Haefliger A., \og Structures feuillet\'ees et cohomologie \`a valeurs dans un faisceau de groupo\"\i des\fg, Comment. Math. Helv., 32, 248-329, \textit{1958}.

\bibitem{Heafliger03} Haefliger A., \og Naissance des feuilletages, d'Ehresmann-Reeb \`a Novikov\fg, pr\'epublication.

\bibitem{HeitschLazarov91} Heitsch J. L., Lazarov C.,  \og Homotopy invariance of foliation Betti numbers\fg,  Invent. Math. 104, no. 2, 321--347, \textit{1991}.


\bibitem{HjorthKechris03} Hjorth G., Kechris A.S., \og Rigidity theorems for actions of product groups and countable Borel equivalence relations\fg, pr\'epublication, \textit{2002}.

\bibitem{JKL} Jackson S., Kechris A.S., Louveau, A. \og Countable borel equivalence relations\fg, J. Math. Log.,  2,  no. 1, 1--80, \textit{2002}.


\bibitem{JonesSchmidt87} Jones V.F.R., Schmidt K., \og Asymptotically invariant sequences and approximate finiteness\fg, Amer. J. Math., 109, 91-114, \textit{1987}.


\bibitem{JuncoRosenblatt79} del Junco A., Rosenblatt J., \og Conterexample in ergodic theory and number theory\fg, Math. Ann., 245, 185-197, \textit{1979}.

\bibitem{Kaima97} Kaimanovich V.A., \og Amenability, hyperfiniteness, and isoperimetric inequalities\fg, C. R. Acad. Sci. Paris S\'er I Math., 325, 999-1004, \textit{1997}.

\bibitem{Ledoux01} Ledoux M., \og The concentration of measure phenomenon\fg,  Mathematical Surveys and Monographs, 89. American Mathematical Society, Providence, RI, \textit{2001}. 


\bibitem{LosertRindler81} Losert V., Rindler H., \og Almost invariant sets\fg,  Bull. London Math. Soc.  13, no. 2, 145--148 \textit{1981}.

\bibitem{Mackey66} Mackey G.W., \og Ergodic theory and virtual groups\fg, Math. Ann., 166, 187-207, \textit{1966}.


\bibitem{Moore82} Moore C.C.,  \og Ergodic theory and von Neumann algebras\fg, In Operator algebras and applications, Part 2 (Kingston, Ont., 1980), 179-226, Amer. Math. Soc., Providence, R.I., \textit{1982}.

\bibitem{OrnsteinWeiss80} Ornstein D., Weiss B., \og Ergodic theory of amenable group action I. The Rohlin lemma \fg, Bull. Amer. Math. Soc., 2(1), 161--164, \textit{1980}.


\bibitem{Pestov00} Pestov V. G., \og Amenable representations and dynamics of the unit sphere in an infinite-dimensional Hilbert space\fg,  Geom. Funct. Anal.,  10,  no. 5, 1171--1201, \textit{2000}.


\bibitem{Pichot04_II} Pichot M., \og Sur les espaces mesur\'es singuliers II\fg, en pr\'eparation.

\bibitem{Pichot05} Pichot M., \og Quasi-p\'eriodicit\'e et th\'eorie de la mesure \fg, Th\`ese, en pr\'eparation.

\bibitem{Ramsay82} Ramsay A., \og Topologies on measured groupoids\fg, J. Funct. Anal., 47, 314-343,  \textit{1982}. 

\bibitem{Rosenblatt81} Rosenblatt J., \og Uniqueness of invariant means for measure preserving transformations\fg,  Trans. Amer. Math. Soc., 265, 623-636, \textit{1981}.


\bibitem{Schmidt80} Schmidt K., \og Asymptotically invariant sequences and an action of $\SL(2,\ZI)$ on the 2-sphere\fg, Israel J. Math., 37, 193-208, \textit{1980}

\bibitem{Schmidt81} Schmidt K., \og Amenability, Kazhdan's property T, strong ergodicity and invariant means for ergodic group action\fg, Erg. Th. and Dyn Systems, 1, 233-236, \textit{1981}.

\bibitem{Zimmer84} Zimmer  R.,  \og Ergodic theory and semi-simple groups\fg, Birkh\"auser Verlag, Basel,\textit{1984}.



\end{thebibliography}
\end{document}